\documentclass[12pt]{article}
\usepackage{amssymb}

\renewcommand{\labelenumi}{(\theenumi)}
\renewcommand{\theenumi}{\roman{enumi}}
\newtheorem{definition}{Definition}[section]
\newtheorem{theorem}[definition]{Theorem}
\newtheorem{lemma}[definition]{Lemma}
\newtheorem{corollary}[definition]{Corollary}

\newtheorem{example}[definition]{Example}

\newtheorem{problem}[definition]{Problem}
\newtheorem{note}[definition]{Note}

\def\R{\mathbb R}
\def\C{\mathbb C}
\def\K{\mathbb K}
\def\D{\mathbb D}
\def\Z{\mathbb Z}
\def\Mdf{{\hbox{Mat}}_{d+1}(\K )}

\begin{document}

\title{ \bf Two linear transformations each tridiagonal \\ 
with respect to an eigenbasis of the other; comments
on the parameter array\footnote{
\hfill\break
\noindent
{\bf 2000 Mathematics Subject Classification}.
05E30, 17B37, 33C45, 33D45. 
}
}
\author{Paul Terwilliger}  

\date{}
\maketitle
\begin{abstract}
Let $\K$ denote a field.
 Let $d$ denote a nonnegative integer and consider a sequence 
$p=(\theta_i, \theta^*_i,i=0...d; \varphi_j, \phi_j,j=1...d)$ consisting
of scalars taken from
$\K$. We call $p$ a {\it parameter array} whenever:
(PA1) 
$\theta_i \not=\theta_j, \;
 \theta^*_i\not=\theta^*_j$ if $i\not=j$,
$(0 \leq i,j\leq d)$;
(PA2)
 $
\varphi_i\not=0$, $\phi_i\not=0$ $(1 \leq i \leq d)$;
(PA3)
$\varphi_i = \phi_1 \sum_{h=0}^{i-1} 
\frac{\theta_h-\theta_{d-h}}{\theta_0-\theta_d}
+ (\theta^*_i-\theta^*_0)(\theta_{i-1}-\theta_d)$
$(1 \leq i \leq d)$;
(PA4)
$\phi_i = \varphi_1 \sum_{h=0}^{i-1} 
\frac{\theta_h-\theta_{d-h}}{\theta_0-\theta_d}
+ (\theta^*_i-\theta^*_0)(\theta_{d-i+1}-\theta_0)$
$(1 \leq i \leq d)$;
(PA5)
$(\theta_{i-2}-\theta_{i+1})(\theta_{i-1}-\theta_i)^{-1}$, 
$(\theta^*_{i-2}-\theta^*_{i+1})(\theta^*_{i-1}-\theta^*_i)^{-1}$
are equal and independent of $i$ for $2 \leq i \leq d-1$.
In \cite{LS99} we showed the 
parameter arrays 
 are in bijection
with the isomorphism classes of Leonard systems. Using
 this bijection
we obtain the following two characterizations of parameter arrays.
Assume $p$ satisfies PA1, PA2. 
Let  $A, B,A^*, B^*$ denote the matrices in 
 $\hbox{Mat}_{d+1}(\K)$ which have  entries
$A_{ii}=\theta_i$,
 $B_{ii}=\theta_{d-i}$,
$A^*_{ii}=\theta^*_i$,
$ B^*_{ii}=\theta^*_{i}$
$(0 \leq i \leq d)$,
$A_{i,i-1}=1$,  
$B_{i,i-1}=1$,
$A^*_{i-1,i}=\varphi_i$, 
$B^*_{i-1,i}=\phi_i$  $(1 \leq i \leq d)$,
and all other entries 0.
We show 
the following are equivalent:
(i) $p$ satisfies PA3--PA5; (ii) there exists an invertible 
$G \in 
\hbox{Mat}_{d+1}(\K)$ such that 
$G^{-1}AG=B$ and 
$G^{-1}A^*G=B^*$;
(iii) for $0 \leq i \leq d$ the polynomial
\begin{eqnarray*}
\sum_{n=0}^i \frac{
(\lambda-\theta_0)
(\lambda-\theta_1) \cdots
(\lambda-\theta_{n-1})
(\theta^*_i-\theta^*_0)
(\theta^*_i-\theta^*_1) \cdots
(\theta^*_i-\theta^*_{n-1})
}
{\varphi_1\varphi_2\cdots \varphi_n}
\end{eqnarray*}
is a scalar multiple of the polynomial 
\begin{eqnarray*}
\sum_{n=0}^i \frac{
(\lambda-\theta_d)
(\lambda-\theta_{d-1}) \cdots
(\lambda-\theta_{d-n+1})
(\theta^*_i-\theta^*_0)
(\theta^*_i-\theta^*_1) \cdots
(\theta^*_i-\theta^*_{n-1})
}
{\phi_1\phi_2\cdots \phi_n}.
\end{eqnarray*}
We display all the parameter arrays in parametric form. For each
array we compute the above polynomials.
The resulting polynomials
form a class consisting of the 
$q$-Racah, $q$-Hahn, dual $q$-Hahn,
$q$-Krawtchouk,
dual $q$-Krawtchouk,
quantum $q$-Krawtchouk, 
affine $q$-Krawtchouk, 
Racah, Hahn, dual-Hahn, Krawtchouk, Bannai/Ito, and  Orphan polynomials.
The Bannai/Ito polynomials can be obtained from the $q$-Racah polynomials by letting $q$ 
tend to $-1$. The Orphan polynomials have maximal degree 3 and 
exist for $\mbox{char}(\K)=2$  only. 
For each of the polynomials listed above
we give the orthogonality,  3-term recurrence, 
and difference equation in terms of the parameter array.
\end{abstract}

\section{Introduction}
In this paper we continue to develop the theory
of Leonard pairs and Leonard systems 
\cite{TD00}, 
\cite{LS99},
\cite{qSerre},
\cite{LS24},
\cite{lsint},
\cite{Terint},
\cite{conform},
\cite{split}.
We briefly summarize our results so far.
In 
\cite{LS99}
we introduced the notion of a Leonard pair and the closely related notion
of a  Leonard system (see Section 2 below.)
We 
classified the Leonard systems. In the process we introduced the
split decomposition for Leonard systems. 
Moreover we showed that
 every Leonard pair satisfies two cubic polynomial relations which
we call the tridiagonal relations.
The tridiagonal relations  generalize both the cubic $q$-Serre relations
and the Dolan-Grady relations.
In 
\cite{TD00} we introduced a generalization of a Leonard pair (resp. system)
which we call a tridiagonal pair (resp. system.)  
We extended  
some of our results on Leonard pairs and systems
to tridiagonal pairs and systems.
For instance we showed 
that every 
tridiagonal system 
 has a split decomposition.
Moreover we showed that
 every tridiagonal pair satisfies
an appropriate pair of tridiagonal relations. 
We did not get a classification of tridiagonal systems and to our knowledge
this remains an open problem.
In 
\cite{qSerre} we introduced the tridiagonal algebra.
This is an associative algebra on two
generators subject to a pair of tridiagonal relations.
 We showed that every tridiagonal pair
induces on the underlying vector space  the structure
of an irreducible  module for a tridiagonal algebra.
Given an irreducible finite dimensional module for a tridiagonal
algebra, we displayed sufficient conditions for it to be induced
from a Leonard pair in this fashion. 
We also showed each sequence of Askey-Wilson polynomials
gives a  basis for an appropriate
infinite dimensional irreducible tridiagonal algebra module.
In \cite{LS24} we began with an arbitrary Leonard pair, and exhibited 24
bases for the underlying vector space which we found attractive. 
For each of these bases we computed the matrices which represent the Leonard pair.
We found each of these
 matrices is tridiagonal, diagonal,
 upper bidiagonal or lower bidiagonal. 
We computed the transition matrix for 
sufficiently many  ordered pairs of bases in our
set of 24 to enable one to readily find the transition matrix for 
any ordered pair of  
bases in our set of 24.
In the survey \cite{lsint} we  gave a number of examples of Leonard pairs. We used these examples
to illustrate how Leonard pairs arise in representation theory, combinatorics,
and the theory of orthogonal polynomials. 
The paper \cite{Terint} is another survey.
In 
\cite{conform}
we introduced the notion of a parameter array. We showed
 that the classification of
 Leonard systems mentioned above  gives 
a bijection from the set of isomorphism
classes of Leonard systems to the set of parameter arrays.
 We  introduced the $TD$-$D$ canonical form and the  $LB$-$UB$ canonical
form  for Leonard  systems.
For a Leonard system in
$TD$-$D$ canonical form the associated Leonard pair is
 represented by a tridiagonal and diagonal
matrix, subject to a certain normalization. For a Leonard system in 
 $LB$-$UB$ canonical form the associated Leonard pair is 
 represented by a lower bidiagonal and upper bidiagonal matrix, subject to
a certain normalization.
We showed every Leonard system is isomorphic to a unique Leonard system which is
in 
$TD$-$D$ canonical form and a unique 
 Leonard system  which is
in 
$LB$-$UB$ canonical form. We described these canonical forms using the 
associated parameter array.
In \cite{split} we obtained two characterizations of Leonard pairs based on
 the split decomposition.

\medskip
\noindent  
We now give an overview of the  present paper.  
We first review our bijection between the set 
of isomorphism classes of Leonard systems
and the set of  parameter arrays.
We then use this bijection to obtain two characterizations
of Leonard systems. 
The first characterization involves bidiagonal matrices
and is given in 
Theorem
\ref{thm:gmat}.
The second characterization involves polynomials 
and is given in 
Theorem 
\ref{eq:mth}.
We view Theorem
\ref{eq:mth} as a variation on a theorem of D. Leonard
\cite[p. 260]{BanIto},
\cite{Leon}.
In Section 5
 we display all the parameter arrays.
For each parameter array we display the corresponding polynomials 
from our second characterization.
These corresponding  polynomials 
form a class consisting of the 
$q$-Racah, $q$-Hahn, dual $q$-Hahn,
$q$-Krawtchouk,
dual $q$-Krawtchouk,
quantum $q$-Krawtchouk, 
affine $q$-Krawtchouk, 
Racah, Hahn, dual-Hahn, Krawtchouk,
 Bannai/Ito, and  Orphan polynomials.
The Bannai/Ito polynomials can be obtained from the $q$-Racah polynomials by letting $q$ 
tend to $-1$. The Orphan polynomials have maximal degree 3 and 
exist for $\mbox{char}(\K)=2$  only. 
For each of the polynomials listed above
we give the orthogonality,  3-term recurrence, 
and difference equation in terms of the parameter array.
We conclude the paper with an open problem.

\medskip
\noindent 
We now recall the definition of a parameter array.
For the rest of this paper $\K$ will denote a field.

\begin{definition} 
\cite[Definition 10.1]{conform}
\label{def:pa}
Let $d$ denote a nonnegative integer.
By a {\it parameter array over $\K$ of diameter $d$}
we mean a sequence of scalars 
$
(\theta_i, \theta^*_i,i=0...d; \varphi_j, \phi_j,j=1...d)
$
taken from $\K$ which satisfy the following conditions
(PA1)--(PA5).
\begin{description}
\item[(PA1)]  
$
\theta_i \not=\theta_j, \qquad 
\theta^*_i \not=\theta^*_j \qquad \mbox{if}\quad i\not=j, \qquad \qquad
(0 \leq i,j\leq d)
$.
\\
\item[(PA2)]
$
\varphi_i \not=0, \qquad \phi_i \not=0 \qquad \qquad (1 \leq i \leq d).
$
\\
\item[(PA3)]
$\varphi_i = \phi_1 \sum_{h=0}^{i-1} 
\frac{\theta_h-\theta_{d-h}}{\theta_0-\theta_d}
+ (\theta^*_i-\theta^*_0)(\theta_{i-1}-\theta_d)
\qquad \qquad (1 \leq i \leq d)$.
\\
\item[(PA4)] 
$\phi_i = \varphi_1 \sum_{h=0}^{i-1} 
\frac{\theta_h-\theta_{d-h}}{\theta_0-\theta_d}
+ (\theta^*_i-\theta^*_0)(\theta_{d-i+1}-\theta_0)
\qquad \qquad (1 \leq i \leq d)$.
\\
\item[(PA5)]
The expressions
\begin{eqnarray} 
\frac{\theta_{i-2}-\theta_{i+1}}{\theta_{i-1}-\theta_i},
\qquad 
\frac{\theta^*_{i-2}-\theta^*_{i+1}}{\theta^*_{i-1}-\theta^*_i}
\label{eq:betaplusone}
\end{eqnarray}
are equal and independent of $i$ for $2 \leq i \leq d-1$.
\end{description}
\end{definition}

\noindent  We now turn our attention to Leonard systems.

\section{Parameter arrays and Leonard systems}

We recall  the notion of a Leonard system and
discuss how these objects are related to parameter arrays.
Our account will be brief; for 
 more detail see \cite{LS99},
\cite{LS24},
\cite{conform},
\cite{lsint}.
Let $d$ denote a nonnegative integer.
Let $\hbox{Mat}_{d+1}(\K)$ denote the $\K$-algebra
consisting of all $d+1$ by $d+1$ matrices which have
entries in $\K$. We index the rows and columns by $0,1,\ldots, d$.
Let $\mathcal A$ denote a $\K$-algebra isomorphic to 
$\hbox{Mat}_{d+1}(\K)$.
An element $A \in {\mathcal A}$ is called {\it multiplicity-free}
whenever it has $d+1$ mutually distinct eigenvalues in
$\K$. Let $A$ denote a multiplicity-free 
element of $\mathcal A$. Let $\theta_0, \theta_1, \ldots, \theta_d$
denote an ordering of the eigenvalues of $A$, and for
$0 \leq i \leq d$ put
\begin{eqnarray*}
E_i = \prod_{{0 \leq j\leq d}\atop {j\not=i}}
\frac{A-\theta_jI}{\theta_i-\theta_j},
\end{eqnarray*}
where $I$ denotes the identity of $\mathcal A$.
We observe $AE_i=\theta_iE_i$ $(0 \leq i \leq d)$;
(ii) $E_iE_j=\delta_{ij}E_i
$ $(0 \leq i,j\leq d)$;
(iii) $\sum_{i=0}^d E_i= I$.
Let $\mathcal D$ denote the subalgebra of $\mathcal A$
generated by $A$. Using (i)--(iii) we find
 $E_0, E_1, \ldots, E_d$ form a basis 
for the $\K$-vector space 
$\mathcal D$.
We call $E_i$ the {\it primitive idempotent} of
$A$ associated with $\theta_i$.
By a {\it Leonard system} in $\cal A$ we mean 
a sequence
$\Phi=(A;A^*; \lbrace E_i\rbrace_{i=0}^d; \lbrace E^*_i\rbrace_{i=0}^d)$
which satisfies the following (i)--(v).
\begin{enumerate}
\item Each of $A,A^*$ is a multiplicity-free element of $\mathcal A$.
\item $E_0, E_1, \ldots, E_d$
is an ordering of the primitive idempotents
of $A$.
\item $E^*_0, E^*_1, \ldots, E^*_d$
is an ordering of the primitive idempotents
of $A^*$.
\item 
${\displaystyle{
E^*_iAE^*_j = \cases{0, &if $\;|i-j|> 1$;\cr
\not=0, &if $\;|i-j|=1$\cr}
\qquad \qquad 
(0 \leq i,j\leq d).
}}$
\item
${\displaystyle{
 E_iA^*E_j = \cases{0, &if $\;|i-j|> 1$;\cr
\not=0, &if $\;|i-j|=1$\cr}
\qquad \qquad 
(0 \leq i,j\leq d).
}}$
\end{enumerate}
We call $\mathcal A$ the {\it ambient algbebra} of
$\Phi$ and  say $\Phi$ is {\it over }$\K$ \cite[Definition 1.4]{LS99}.

\medskip
\noindent 
Let 
$\Phi=(A;A^*; \lbrace E_i\rbrace_{i=0}^d; \lbrace E^*_i\rbrace_{i=0}^d)$
denote a Leonard system in $\mathcal A$. 
Then each of the following is a Leonard system in $\mathcal A$:
\begin{eqnarray*}
\Phi^*&:=&(A^*;A; \lbrace E^*_i\rbrace_{i=0}^d; \lbrace E_i\rbrace_{i=0}^d),
\\
\Phi^\downarrow &:=&(A;A^*; \lbrace E_i\rbrace_{i=0}^d; \lbrace E^*_{d-i}\rbrace_{i=0}^d),
\\
\Phi^\Downarrow &:=&(A;A^*; \lbrace E_{d-i}\rbrace_{i=0}^d; \lbrace E^*_{i}\rbrace_{i=0}^d).
\end{eqnarray*}
Viewing $*, \downarrow, \Downarrow$ as permutations on the set of all
Leonard systems,
\begin{eqnarray}
&&\qquad \qquad *^2=\downarrow^2 = \Downarrow^2 = 1, 
\label{eq:d4rel1}
\\
&&\Downarrow *= * \downarrow, \qquad \qquad\downarrow *= *\Downarrow, 
\qquad \qquad \downarrow \Downarrow = \Downarrow \downarrow.
\label{eq:d4rel2}
\end{eqnarray}
The group generated by the symbols  
 $*, \downarrow, \Downarrow$ subject to the relations
(\ref{eq:d4rel1}), (\ref{eq:d4rel2}) is the dihedral group $D_4$.
We recall $D_4$ is the group of symmetries of a square, and has 
8 elements. Apparently 
 $*, \downarrow, \Downarrow$ induce an action of $D_4$ on the set
 of all Leonard systems.

\medskip
\noindent 
Let 
$\Phi=(A;A^*; \lbrace E_i\rbrace_{i=0}^d; \lbrace E^*_i\rbrace_{i=0}^d)$
denote a Leonard system in $\mathcal A$. 
In order to describe $\Phi$ we define some
parameters.
For $0 \leq i \leq d$ let $\theta_i$ (resp. $\theta^*_i$)
denote the eigenvalue of $A$ (resp. $A^*$)
associated with $E_i$ (resp. $E^*_i$.) We call
$\theta_0, \theta_1, \ldots, \theta_d$
(resp. 
$\theta^*_0, \theta^*_1, \ldots, \theta^*_d$)
the {\it eigenvalue sequence} (resp. 
{\it dual eigenvalue sequence})
of $\Phi$.
Let $V$ denote an irreducible left $\mathcal A$-module.
By a decomposition of $V$
we mean a sequence $U_0, U_1, \ldots, U_d$
 consisting of 1-dimensional
subspaces of $V$ such that
\begin{eqnarray*}
V= U_0+U_1+\cdots +U_d \qquad \qquad  \mbox{(direct sum)}.
\end{eqnarray*}
By \cite[Theorem 3.2]{LS99} there exists a unique decomposition 
$U_0, U_1, \ldots, U_d$ of $V$ such that both
\begin{eqnarray}
(A-\theta_iI)U_i&=& U_{i+1} \qquad (0 \leq i \leq d-1), 
\quad (A-\theta_dI)U_d=0,
\label{eq:split1}
\\
(A^*-\theta^*_iI)U_i&=& U_{i-1} \qquad (1 \leq i \leq d), 
\quad (A^*-\theta^*_0I)U_0=0.
\label{eq:split2}
\end{eqnarray}
Pick any integer $i$ $(1 \leq i\leq d)$. Then 
$(A^*-\theta^*_iI)U_i=U_{i-1}$ and 
$(A-\theta_{i-1}I)U_{i-1}=U_{i}$.
Apparently $U_i$ is an eigenspace for 
$(A-\theta_{i-1}I)
(A^*-\theta^*_iI)$ and the corresponding eigenvalue is a nonzero
scalar in $\K$. We denote this eigenvalue by $\varphi_i$.
We call  $\varphi_1, \varphi_2, \ldots, \varphi_d$ the 
{\it first split sequence} of 
$
\Phi$.
We let $\phi_1, \phi_2, \ldots, \phi_d$ denote the first split
sequence of $\Phi^\Downarrow$
and call this the {\it second split sequence} of 
$\Phi$.

\medskip
\noindent We recall the notion of {\it isomorphism}
for Leonard systems.  
Let 
$\Phi=(A;A^*; \lbrace E_i\rbrace_{i=0}^d; \lbrace E^*_i\rbrace_{i=0}^d)$
denote a Leonard system in $\mathcal A$ and let 
$\sigma : {\mathcal A}\rightarrow {\mathcal A}'$ denote an isomorphism
of $\K$-algebras. We write
$\Phi^\sigma=(A^\sigma;A^{*\sigma}; \lbrace E^{\sigma}_i\rbrace_{i=0}^d; \lbrace E^{*\sigma}_i\rbrace_{i=0}^d)$
and observe $\Phi^{\sigma}$ is a Leonard system in ${\mathcal A}'$.
Let $\Phi$
and 
$\Phi'$
denote any Leonard systems over $\K$.
By an {\it isomorphism of Leonard systems} from $\Phi$ to $\Phi'$ we mean
an isomorphism of $\K$-algbras from the ambient algebra of $\Phi$
to the ambient algebra of $\Phi'$ such that
${\Phi}^\sigma = \Phi'$.
We say $\Phi$ and $\Phi'$ are {\it isomorphic} whenever
there exists an isomorphism of Leonard systems from
$\Phi$ to $\Phi'$.

\begin{theorem} \cite[Theorem 1.9]{LS99}
\label{thm:ls}
Let $d$ denote a nonnegative  integer and let 
$(\theta_i, \theta^*_i,i=0...d;\varphi_j,\phi_j,j=1...d)$
denote a sequence of scalars taken from $\K$.
Then the following (i), (ii) are equivalent.
\begin{enumerate}
\item The sequence 
$(\theta_i, \theta^*_i,i=0...d;\varphi_j,\phi_j,j=1...d)$
is a parameter array over $\K$.
\item
There exists 
a Leonard system $\Phi$ over $\K$ which has
eigenvalue sequence 
$\theta_0, \theta_1, \ldots, \theta_d$, dual eigenvalue sequence
$\theta^*_0, \theta^*_1, \ldots, \theta^*_d$,
first split sequence 
$\varphi_1, \varphi_2, \ldots, \varphi_d$ and second split sequence
$\phi_1, \phi_2, $ $ \ldots, \phi_d$.
\end{enumerate}
Suppose (i), (ii) hold. Then $\Phi$ is 
unique up to isomorphism of Leonard systems.
\end{theorem}

\noindent 
Let 
$\Phi=(A;A^*; \lbrace E_i\rbrace_{i=0}^d; \lbrace E^*_i\rbrace_{i=0}^d)$
denote a Leonard system. By the {\it parameter array of $\Phi$}
we mean the sequence
$(\theta_i, \theta^*_i,i=0...d;\varphi_j,\phi_j,j=1...d)$
where 
$\theta_0, \theta_1, \ldots, \theta_d$
(resp. 
$\theta^*_0, \theta^*_1, \ldots, \theta^*_d$)
is the eigenvalue sequence
(resp. 
dual  eigenvalue sequence)
of $\Phi$ and
$\varphi_1, \varphi_2, \ldots, \varphi_d$
(resp.
$\phi_1, \phi_2, $ $ \ldots, \phi_d$)
is the first split sequence 
(resp.
 second split sequence )
of $\Phi$.
By Theorem \ref{thm:ls} the map which sends a given Leonard system
to its parameter array induces a bijection from the set of isomorphism
classes of Leonard systems over $\K$ to the set of parameter arrays over
$\K$. 

\medskip
\noindent Earlier we mentioned an action of $D_4$ on the set of 
Leonard systems. The above bijection induces an action of $D_4$
on the set of parameter arrays. This action is described as follows.

\begin{lemma}\cite[Theorem 1.11]{LS99}
\label{lem:d4}
Let
$\Phi=(A;A^*; \lbrace E_i\rbrace_{i=0}^d; \lbrace E^*_i\rbrace_{i=0}^d)$
denote a Leonard system
and let  $p=
(\theta_i, \theta^*_i,i=0...d;\varphi_j,\phi_j,j=1...d)$
denote the corresponding parameter array.
\begin{enumerate}
\item The parameter array of $\Phi^*$ is $p^*$ where
$p^*:=(\theta^*_i, \theta_i,i=0...d;\varphi_j,\phi_{d-j+1},j=1...d)$.
\item The parameter array of $\Phi^\downarrow$ is $p^\downarrow$ where
$p^\downarrow:=(\theta_i, \theta^*_{d-i},i=0...d;\phi_{d-j+1},\varphi_{d-j+1},j=1...d)$.
\item The parameter array of $\Phi^\Downarrow$ is $p^\Downarrow$ where
$p^\Downarrow:=(\theta_{d-i}, \theta^*_i,i=0...d;\phi_{j},\varphi_{j},j=1...d)$.
\end{enumerate}
\end{lemma}

\section{Parameter arrays and bidiagonal matrices}

In this section we characterize the parameter
arrays in terms of bidiagonal matrices.
We will refer to the following set-up.

\begin{definition}
\label{def:setup}
Let $d$ denote a nonnegative integer and let
$(\theta_i, \theta^*_i,i=0...d; \varphi_j, \phi_j,j=1...d)
$
denote a sequence of  scalars taken from $\K$. We assume
this sequence satisfies 
PA1 and PA2.
\end{definition}

\begin{theorem}
\label{thm:gmat}
With reference to Definition \ref{def:setup},
 the following (i), (ii) are equivalent.
\begin{enumerate}
\item
The sequence
$(\theta_i, \theta^*_i,i=0...d;\varphi_j,\phi_j,j=1...d)$
satisfies PA3-PA5.
\item
There exists an invertible matrix 
$G \in \mbox{Mat}_{d+1}(\K)$ such that both
\begin{eqnarray}
G^{-1}
\left(
\begin{array}{c c c c c c}
\theta_0 & & & & & {\mathbf 0} \\
1 & \theta_{1} &  & & & \\
& 1 & \theta_{2} &  & & \\
& & \cdot & \cdot &  &  \\
& & & \cdot & \cdot &  \\
{\mathbf 0}& & & & 1 & \theta_d
\end{array}
\right) G
&=&\left(
\begin{array}{c c c c c c}
\theta_d & & & & & {\mathbf 0} \\
1 & \theta_{d-1} &  & & & \\
& 1 & \theta_{d-2} &  & & \\
& & \cdot & \cdot &  &  \\
& & & \cdot & \cdot &  \\
{\mathbf 0}& & & & 1 & \theta_0
\end{array}
\right),
\label{eq:g1}
\\
G^{-1}
\left(
\begin{array}{c c c c c c}
\theta^*_0 &\varphi_1 & & & & {\mathbf 0} \\
 & \theta^*_1 & \varphi_2 & & & \\
&  & \theta^*_2 & \cdot & & \\
& &  & \cdot & \cdot &  \\
& & &  & \cdot & \varphi_d \\
{\mathbf 0}& & & &  & \theta^*_d
\end{array}
\right) G
&=&
\left(
\begin{array}{c c c c c c}
\theta^*_0 &\phi_1 & & & & {\mathbf 0} \\
 & \theta^*_1 & \phi_2 & & & \\
&  & \theta^*_2 & \cdot & & \\
& &  & \cdot & \cdot &  \\
& & &  & \cdot & \phi_d \\
{\mathbf 0}& & & &  & \theta^*_d
\end{array}
\right).
\label{eq:g2}
\end{eqnarray}
\end{enumerate}
\end{theorem}
\noindent {\it Proof:}
$(i)\Rightarrow (ii)$
The sequence
$(\theta_i, \theta^*_i,i=0...d;\varphi_j,\phi_j,j=1...d)$
satisfies PA1--PA5 and is therefore a parameter array over $\K$.  
By Theorem 
\ref{thm:ls} there exists a Leonard system
over $\K$ which has
 eigenvalue sequence 
$\theta_0, \theta_1, \ldots, \theta_d$, dual eigenvalue sequence
$\theta^*_0, \theta^*_1, \ldots, \theta^*_d$,
first split sequence 
$\varphi_1, \varphi_2, \ldots, \varphi_d$ and second split sequence
$\phi_1, \phi_2, $ $ \ldots, \phi_d$.
We denote this system by 
$\Phi=(A;A^*; \lbrace E_i \rbrace_{i=0}^d; \lbrace E^*_i\rbrace_{i=0}^d )$.
Let $\mathcal A$ denote the ambient algebra of
$\Phi$ and
let $V$ denote an irreducible left $\mathcal A$-module.
Let $U_0, U_1, \ldots, U_d$ denote the decomposition of $V$
which satisfies
(\ref{eq:split1}), (\ref{eq:split2}).
For $0 \leq i \leq d$ let $u_i$ denote a nonzero vector in
$U_i$ and observe $u_0, u_1, \ldots, u_d$ is a basis for
$V$. Normalizing this basis we may assume
 $(A-\theta_iI)u_i = u_{i+1}$ for $0 \leq i \leq d-1$ and
 $(A-\theta_dI)u_d =0$. Since 
$\varphi_1, \varphi_2, \ldots, \varphi_d$ is the first split
sequence of $\Phi$ we have
$(A^*-\theta^*_iI)u_i = \varphi_iu_{i-1}$ $(1 \leq i \leq d)$,
 $(A^*-\theta^*_0I)u_0 =0$.
Applying these comments to $\Phi^\Downarrow$
we find there exists 
 a basis $v_0, v_1, \ldots, v_d$
for $V$ such that $(A-\theta_{d-i})v_i = v_{i+1}$ $(0 \leq i \leq d-1)$,
 $(A-\theta_0)v_d =0$ and
$(A^*-\theta^*_i)v_i = \phi_iv_{i-1}$ $(1 \leq i \leq d)$,
 $(A^*-\theta^*_0)v_0 =0$. Let $G \in \mbox{Mat}_{d+1}(\K)$
 denote the transition matrix
 from the basis  
 $u_0, u_1, \ldots, u_d$
to the  
 basis $v_0, v_1, \ldots, v_d$, so that
$
v_j = \sum_{i=0}^d G_{ij}u_i$ for $0 \leq j \leq d$.
Using this and elementary linear algebra 
we find $G$ is invertible and satisfies 
(\ref{eq:g1}), (\ref{eq:g2}).
\\
$(ii)\Rightarrow (i)$
We apply Theorem \ref{thm:ls}. We show
condition (ii) holds in that theorem.
In order to do this we invoke some results from
\cite{split}.
Consider the following matrices in 
$\mbox{Mat}_{d+1}(\K)$:
$$
A=\left(
\begin{array}{c c c c c c}
\theta_0 & & & & & {\mathbf 0} \\
1 & \theta_{1} &  & & & \\
& 1 & \theta_{2} &  & & \\
& & \cdot & \cdot &  &  \\
& & & \cdot & \cdot &  \\
{\mathbf 0}& & & & 1 & \theta_d
\end{array}
\right) 
,
\qquad 
A^*=
\left(
\begin{array}{c c c c c c}
\theta^*_0 &\varphi_1 & & & & {\mathbf 0} \\
 & \theta^*_1 & \varphi_2 & & & \\
&  & \theta^*_2 & \cdot & & \\
& &  & \cdot & \cdot &  \\
& & &  & \cdot & \varphi_d \\
{\mathbf 0}& & & &  & \theta^*_d
\end{array}
\right).
$$
We observe $A$ (resp. $A^*$) is multiplicity-free,
with eigenvalues $\theta_0, \theta_1,\ldots, \theta_d$
(resp. $\theta^*_0, \theta^*_1,\ldots,$ $ \theta^*_d$.)
For $0 \leq i \leq d$ we let $E_i$  (resp. $E^*_i$) denote
the primitive idempotent for $A$ (resp. $A^*$) associated
with $\theta_i$ (resp. $\theta^*_i$.)
By \cite[Lemma 6.2, Theorem 6.3]{split} 
 the sequence
$(A;A^*; \lbrace E_i \rbrace_{i=0}^d; $ $\lbrace E^*_i\rbrace_{i=0}^d )$  
is a Leonard system 
in $\mbox{Mat}_{d+1}(\K)$. 
Let us call this system $\Phi$. 
By the construction 
$\Phi$ has eigenvalue sequence 
 $\theta_0, \theta_1,\ldots, \theta_d$
and dual eigenvalue sequence $\theta^*_0, \theta^*_1,\ldots,\theta^*_d$.
From the form of $A$ and $A^*$ we find
$\varphi_1, \varphi_2, \ldots, \varphi_d$
is the first split sequence for $\Phi$. By the last line of
\cite[Theorem 6.3]{split} we find 
$\phi_1, \phi_2, \ldots, \phi_d$ is the second split sequence for
$\Phi$. Now Theorem \ref{thm:ls}(ii) holds; applying
that theorem  we find
 $(\theta_i, \theta^*_i,i=0...d; \varphi_j, \phi_j,j=1...d)$
is a parameter array over $\K$. In particular 
 $(\theta_i, \theta^*_i,i=0...d; \varphi_j, \phi_j,j=1...d)$
satisfies PA3--PA5.
\hfill $\Box $ \\

\noindent The matrix $G$ from Theorem 
\ref{thm:gmat}(ii) will be discussed further in Section
10.

\section{Parameter arrays and polynomials}

\noindent In this section we characterize the parameter
arrays in terms of polynomials.
We will use the following notation.
Let $\lambda $ denote an indeterminate, and let
$\K \lbrack \lambda \rbrack $ denote the $\K$-algebra
consisting of all polynomials in $\lambda $  which have
coefficients in $\K$. For the rest of this paper 
all polynomials which we discuss
are assumed to lie in $\K\lbrack \lambda \rbrack $. 

We view the following theorem as a variation on a theorem
of 
D. Leonard
\cite[p. 260]{BanIto},
\cite{Leon}.

\begin{theorem}
\label{eq:mth}
With reference to Definition \ref{def:setup},
the following (i), (ii) are equivalent.
\begin{enumerate}
\item
The sequence
 $(\theta_i, \theta^*_i,i=0...d; \varphi_j, \phi_j,j=1...d)$
satisfies PA3--PA5.
\item For $0 \leq i \leq d$ the 
polynomial 
\begin{equation}
\sum_{n=0}^i \frac{
(\lambda-\theta_0)
(\lambda-\theta_1) \cdots
(\lambda-\theta_{n-1})
(\theta^*_i-\theta^*_0)
(\theta^*_i-\theta^*_1) \cdots
(\theta^*_i-\theta^*_{n-1})
}
{\varphi_1\varphi_2\cdots \varphi_n}
\label{eq:poly1}
\end{equation}
is a scalar multiple of the polynomial 
\begin{equation}
\sum_{n=0}^i \frac{
(\lambda-\theta_d)
(\lambda-\theta_{d-1}) \cdots
(\lambda-\theta_{d-n+1})
(\theta^*_i-\theta^*_0)
(\theta^*_i-\theta^*_1) \cdots
(\theta^*_i-\theta^*_{n-1})
}
{\phi_1\phi_2\cdots \phi_n}.
\label{eq:poly2}
\end{equation}
\end{enumerate}
\end{theorem}
\noindent {\it Proof:}
Let us abbreviate
$$
A=\left(
\begin{array}{c c c c c c}
\theta_0 & & & & & {\mathbf 0} \\
1 & \theta_{1} &  & & & \\
& 1 & \theta_{2} &  & & \\
& & \cdot & \cdot &  &  \\
& & & \cdot & \cdot &  \\
{\mathbf 0}& & & & 1 & \theta_d
\end{array}
\right) 
,
\qquad 
A^*=
\left(
\begin{array}{c c c c c c}
\theta^*_0 &\varphi_1 & & & & {\mathbf 0} \\
 & \theta^*_1 & \varphi_2 & & & \\
&  & \theta^*_2 & \cdot & & \\
& &  & \cdot & \cdot &  \\
& & &  & \cdot & \varphi_d \\
{\mathbf 0}& & & &  & \theta^*_d
\end{array}
\right)
$$
and 
$$
B=\left(
\begin{array}{c c c c c c}
\theta_d & & & & & {\mathbf 0} \\
1 & \theta_{d-1} &  & & & \\
& 1 & \theta_{d-2} &  & & \\
& & \cdot & \cdot &  &  \\
& & & \cdot & \cdot &  \\
{\mathbf 0}& & & & 1 & \theta_0
\end{array}
\right) 
,
\qquad 
B^*=
\left(
\begin{array}{c c c c c c}
\theta^*_0 &\phi_1 & & & & {\mathbf 0} \\
 & \theta^*_1 & \phi_2 & & & \\
&  & \theta^*_2 & \cdot & & \\
& &  & \cdot & \cdot &  \\
& & &  & \cdot & \phi_d \\
{\mathbf 0}& & & &  & \theta^*_d
\end{array}
\right).
$$
We let
$T, T^*, T^{\Downarrow}$ denote the matrices in
$\mbox{Mat}_{d+1}(\K)$ which have entries
 $T_{ij}= \prod_{h=0}^{j-1}(\theta_i-\theta_h)$,
$T^*_{ij}= \prod_{h=0}^{j-1}(\theta^*_i-\theta^*_h)$,
$T^\Downarrow_{ij}= \prod_{h=0}^{j-1}(\theta_{d-i}-\theta_{d-h})$
for $0 \leq i,j\leq d$.
Each of  
$T, T^*, T^{\Downarrow}$ 
 is lower triangular
 with diagonal entries nonzero so these matrices
 are 
invertible.
Let 
 $D$ (resp. $D^{\Downarrow}$)  denote the diagonal matrix
 in 
$\mbox{Mat}_{d+1}(\K)$ which has ${ii}$th entry
$\varphi_1\varphi_2\cdots \varphi_i$ (resp. 
$\phi_1\phi_2\cdots \phi_i$) for $0 \leq i \leq d$.
Each of $D$, $D^{\Downarrow}$ is invertible.
We let $Z$ denote the matrix in 
$\mbox{Mat}_{d+1}(\K)$ which has ${ij}$th entry
$1$ if $i+j=d$ and 0 if $i+j\not=d$, for $0 \leq i,j\leq d$. 
Observe $Z^2=I$ so 
  $Z$ is invertible.
We let $H$ (resp. $H^*$) denote the diagonal matrix in
$\mbox{Mat}_{d+1}(\K)$ which has ${ii}$th entry
$\theta_i$ 
(resp. $\theta^*_i$) for $0 \leq i \leq d$. 
One verifies $TA=HT$ so $A=T^{-1}HT$.
One verifies 
$ZT^{\Downarrow}B=HZT^{\Downarrow}$ 
so $B=T^{\Downarrow -1}ZHZT^{\Downarrow}$.
One verifies $DA^*D^{-1}T^{*t}=T^{*t}H^*$ so 
$A^*=D^{-1}T^{*t}H^*T^{*-1t}D$.
Similarly 
$B^*=D^{\Downarrow -1}T^{*t}H^*T^{*-1t}D^{\Downarrow}$.
For $0 \leq i \leq d$ let 
$f_i$
denote the polynomial
in (\ref{eq:poly1}).
Let ${\mathcal P}$ denote the matrix in $\mbox{Mat}_{d+1}(\K)$
which has $ij$th entry $f_j(\theta_i)$ for $0 \leq i,j\leq d$.
From the form of
 (\ref{eq:poly1}) 
we find
${\mathcal P}= TD^{-1}T^{*t}$.
For $0 \leq i\leq d$ 
let $f^{\Downarrow}_i$
denote the polynomial
in (\ref{eq:poly2}).
Let ${\mathcal P}^{\Downarrow}$ denote the matrix in 
$\mbox{Mat}_{d+1}(\K)$ which has
 $ij$th entry $f^{\Downarrow}_j(\theta_i)$ for $0 \leq i,j\leq d$.
From the form of
 (\ref{eq:poly2}) 
we find
${\mathcal P}^{\Downarrow}= ZT^{\Downarrow}D^{\Downarrow -1}T^{*t}$.
\\
$(i)\Rightarrow (ii)$
By Theorem
\ref{thm:gmat} there exists an invertible
matrix $G \in 
 \mbox{Mat}_{d+1}(\K)$
such that $G^{-1}AG=B$ and $G^{-1}A^*G=B^*$.
Evaluating
$G^{-1}A^*G=B^*$ using
$A^*=D^{-1}T^{*t}H^*T^{*-1t}D$
and $B^*=D^{\Downarrow -1}T^{*t}H^*T^{*-1t}D^{\Downarrow}$
we find $T^{*-1t}DGD^{\Downarrow -1}T^{*t}$ commutes with $H^*$.
Since $H^*$ is diagonal with diagonal
entries mutually distinct we find
 $T^{*-1t}DGD^{\Downarrow -1}T^{*t}$ is diagonal.
We denote this diagonal matrix by $F$ and
observe
 $G=
D^{-1}T^{*t}FT^{*-1t}D^{\Downarrow}$. 
In this product each factor is upper triangular (or diagonal)
so $G$ is upper triangular.
Recall 
$G^{-1}AG=B$; evaluating this 
using
$A=T^{-1}HT$
and $B=T^{\Downarrow -1}ZHZT^{\Downarrow}$
we find $TGT^{\Downarrow -1}Z$ commutes with
$H$. Since $H$ is diagonal with diagonal
entries mutually distinct we find
 $TGT^{\Downarrow -1}Z$ is diagonal.
We denote this diagonal matrix by $Y$ and observe
$TG=YZT^{\Downarrow}$. In this equation 
we compute the entries in 
column 0.
To aid in this calculation we 
recall $G$ is upper triangular and observe
 $T_{i0}=1$, $T^{\Downarrow}_{i0}=1$ for
$0 \leq i \leq d$. Computing the column 0 entries
in $TG=YZT^{\Downarrow}$ using these facts we find
$Y_{ii}=G_{00}$ for $0 \leq i \leq d$.
Apparently $Y=G_{00}I$.
We remark $G_{00}\not=0$ since $Y$ is invertible by
the construction.
Dividing $G$ by $G_{00}$ we may assume
$G_{00}=1$. Now $Y=I$ so
$G=T^{-1}ZT^{\Downarrow}$.
Recall $T^{*-1t}DGD^{\Downarrow -1}T^{*t}$ is diagonal.
We evaluate this expression using 
$G=T^{-1}ZT^{\Downarrow}$ and find
$T^{*-1t}D  
T^{-1}ZT^{\Downarrow}
D^{\Downarrow -1}T^{*t}$ is diagonal.
But 
$T^{*-1t}D  
T^{-1}ZT^{\Downarrow}
D^{\Downarrow -1}T^{*t}
={\mathcal P}^{-1}{\mathcal P}^\Downarrow $
so 
 ${\mathcal P}^{-1}{\mathcal P}^\Downarrow$ is diagonal.
Taking the inverse 
 we find ${\mathcal P}^{\Downarrow -1}{\mathcal P}$ is diagonal.
For $0 \leq i \leq d$ let $\alpha_i$ denote the $ii$ entry of
this diagonal matrix.
From the definition of
 ${\mathcal P}$ and 
 ${\mathcal P}^\Downarrow$
we find $f_i(\theta_j)= \alpha_i f_i^{\Downarrow}(\theta_j)$
for $0 \leq i ,j\leq d$.
Recall $\theta_0, \theta_1, \ldots, \theta_d$ are mutually
distinct, and that
each of $f_i, f^\Downarrow_i$ has degree $i$ for $0 \leq i \leq d$.
From these comments we find
$f_i = \alpha_if^{\Downarrow}_i$ for $0 \leq i \leq d$.
\\
$(ii)\Rightarrow (i)$
We show Theorem
\ref{thm:gmat}(ii) holds. To do this
we exhibit an invertible matrix 
$G 
\in \mbox{Mat}_{d+1}(\K)$
such that $AG=GB$ and 
 $A^*G=GB^*$.
We define $G=T^{-1}ZT^{\Downarrow}$.
Observe $G$ is invertible.
The equation $AG=GB$ is routinely verified by
evaluating $A,B, G$ using  
 $A=T^{-1}HT$,
$B=T^{\Downarrow -1}ZHZT^{\Downarrow}$,
 $G=T^{-1}ZT^{\Downarrow}$.
We now show $A^*G=GB^*$.
For $0 \leq i \leq d$ there exists
$\alpha_i \in \K$ such that
$f_i=\alpha_i f^{\Downarrow}_i$.
It follows
${\mathcal P}
={\mathcal P}^{\Downarrow}
\mbox{diag}(\alpha_0, \alpha_1, \ldots, \alpha_d)$.
From this and since $H^*$ is diagonal we find
${\mathcal P}
H^*
{\mathcal P}^{-1}
=
{\mathcal P}^{\Downarrow}
H^*
{\mathcal P}^{\Downarrow -1}
$.
In this equation we multiply both sides on the left by $T^{-1}$ 
and on the right by $ZT^{\Downarrow}$
to obtain
$T^{-1}
{\mathcal P}
H^*
{\mathcal P}^{-1}
ZT^{\Downarrow}
=
T^{-1}
{\mathcal P}^{\Downarrow}
H^*
{\mathcal P}^{\Downarrow -1}
ZT^{\Downarrow}
$.
In this equation the left side is equal to 
$A^*G$ and the right side is equal to $GB^*$
so 
$A^*G=GB^*$.
We have now shown $G$ satisfies Theorem
\ref{thm:gmat}(ii). Applying that theorem
we find
 $(\theta_i, \theta^*_i,i=0...d; \varphi_j, \phi_j,j=1...d)$
satisifes PA3--PA5.
\hfill $\Box $ \\

\noindent We finish this section with a comment.

\begin{lemma}
\label{lem:scal}
Referring to Theorem 
\ref{eq:mth}, assume the equivalent conditions (i), (ii)
from that theorem hold. Then for $ 0 \leq i \leq d$ the
scalar referred to in condition (ii) is equal to
\begin{eqnarray*}
\frac{\phi_1\phi_2\cdots \phi_i}{
\varphi_1\varphi_2\cdots \varphi_i}.
\end{eqnarray*}
\end{lemma}
\noindent {\it Proof:}
Compare the coefficient of $\lambda^i$ in
(\ref{eq:poly1}),
(\ref{eq:poly2}).
\hfill $\Box $

\section{The parameter arrays}

In this section we display all the parameter arrays over $\K$.
We will use the following notatation.

\begin{definition}
\label{def:base}
Let
$p=(\theta_i, \theta^*_i, i=0...d;  \varphi_j, \phi_j, j=1...d)$
denote a parameter array over $\K$.  By a {\it base} for
$p$, we mean a nonzero scalar $q$ in the algebraic closure
of $\K$  such that
$q+q^{-1}+1$ is equal to the common value of 
(\ref{eq:betaplusone}) for $2 \leq i \leq d-1$.
We remark on the uniqueness of the base.
Suppose $d\geq 3$.
If $q$ is a base for $p$ then so is $q^{-1}$ and 
$p$ has no other base. Suppose $d<3$.
Then any nonzero scalar in the algebraic closure of 
$\K$ is a base for $p$.
\end{definition}

\begin{definition}
\label{def:assocpoly}
Let
$p=(\theta_i, \theta^*_i, i=0...d;  \varphi_j, \phi_j, j=1...d)$
denote a parameter array over $\K$. 
For $0 \leq i\leq d$ we let $f_i$ denote the following polynomial
in
$\K\lbrack \lambda \rbrack$.
\begin{eqnarray}
f_i = \sum_{n=0}^i \frac{
(\lambda-\theta_0)
(\lambda-\theta_1) \cdots
(\lambda-\theta_{n-1})
(\theta^*_i-\theta^*_0)
(\theta^*_i-\theta^*_1) \cdots
(\theta^*_i-\theta^*_{n-1})
}
{\varphi_1\varphi_2\cdots \varphi_n}.
\label{eq:fipoly}
\end{eqnarray}
We call $f_0, f_1, \ldots, f_d$ the polynomials which {\it correspond } to $p$.
\end{definition}
\noindent 
We now display all the  parameter arrays over $\K$.
For each displayed array 
$(\theta_i, \theta^*_i, i=0...d;  \varphi_j, \phi_j, j=1...d)$
we give a base
and present  $f_i(\theta_j)$
for $0 \leq i,j\leq d$, where  $f_0, f_1, \ldots, f_d$ are the corresponding
polynomials.
Our  
presentation is organized as follows.
In each of Example 
\ref{ex:pa2}--\ref{ex:orphan} below we give a family of
parameter arrays over $\K$. 
In Theorem
\ref{thm:thatsit}
 we show every parameter array over $\K$ is contained in
 at least one of these families.

\medskip
\noindent 
In each of
 Example 
\ref{ex:pa2}--\ref{ex:orphan} below 
the following implicit assumptions apply:
$d$ denotes a nonnegative integer,
the scalars
$(\theta_i, \theta^*_i, i=0...d;  \varphi_j, \phi_j, j=1...d)$
are contained in  $\K$, and the scalars $q,h,h^*\ldots $
are contained in the algebraic closure of $\K$.

\begin{example}($q$-Racah)
\label{ex:pa2}
\label{ex:qracah}
Assume
\begin{eqnarray}
\theta_i &=& \theta_0+h(1-q^i)(1-sq^{i+1})q^{-i},
\label{eq:thrac}
\\
\theta^*_i &=& \theta^*_0+h^*(1-q^i)(1-s^*q^{i+1})q^{-i}
\label{eq:thsrac}
\end{eqnarray}
for $0 \leq i \leq d$ and
\begin{eqnarray}
\varphi_i &=& hh^*q^{1-2i}(1-q^i)(1-q^{i-d-1})(1-r_1q^i)(1-r_2q^i),
\label{lem:varphiphi1}
\\
\phi_i &=& hh^*q^{1-2i}(1-q^i)(1-q^{i-d-1})(r_1-s^*q^i)(r_2-s^*q^i)/s^*
\label{lem:varphiphi2}
\end{eqnarray}
for $1 \leq i \leq d$.
Assume 
$h, h^*, q, s, s^*, r_1, r_2$ are nonzero
 and  $r_1r_2=ss^*q^{d+1}$.
Assume
none of
$q^i,
r_1q^i, r_2q^i,
s^*q^i/r_1,$ $ s^*q^i/r_2$
is equal to $1$ for $1 \leq i \leq d$ and 
that neither of $sq^i, s^*q^i$ is equal to $1$ for $2 \leq i \leq 2d$.
Then 
$(\theta_i, \theta^*_i, i=0...d;  \varphi_j, \phi_j, j=1...d)$
is a parameter array over $\K$ which has base $q$.
The corresponding polynomials
$f_i$  satisfy
\begin{eqnarray*}
f_i(\theta_j)=
{}_4\phi_3 \Biggl({{q^{-i}, \;s^*q^{i+1},\;q^{-j},\;sq^{j+1}}\atop
{r_1q,\;\;r_2q,\;\;q^{-d}}}\;\Bigg\vert \; q,\;q\Biggr)
\label{eq:fqrac}
\end{eqnarray*}
for $0 \leq i,j\leq d$. These $f_i$ are $q$-Racah polynomials.
\end{example}

\begin{example}
($q$-Hahn)
\label{ex:qhahn}
Assume
\begin{eqnarray*}
\theta_i &=&\theta_0+h(1-q^i)q^{-i}, 
\\
\theta^*_i &=& \theta^*_0+h^*(1-q^i)(1-s^*q^{i+1})q^{-i}
\end{eqnarray*}
for $0 \leq i \leq d$ and
\begin{eqnarray*}
\varphi_i &=& hh^*q^{1-2i}(1-q^i)(1-q^{i-d-1})(1-rq^i),
\\
\phi_i &=& -hh^*q^{1-i}(1-q^i)(1-q^{i-d-1})
(r-s^*q^i)
\end{eqnarray*}
for $1 \leq i \leq d$.
Assume 
$h, h^*, q, s^*, r$ are nonzero.
Assume none of
$q^i,
rq^i,
s^*q^i/r$
is equal to $1$ for $1 \leq i \leq d$ and 
that $s^*q^i\not=1$ for $2 \leq i \leq 2d$.
Then the sequence
$(\theta_i, \theta^*_i, i=0...d;  \varphi_j, \phi_j, j=1...d)$
is a parameter array over $\K$ which has base $q$.
The corresponding polynomials
$f_i$
satisfy
\begin{eqnarray*}
f_i(\theta_j)=
{}_3\phi_2 \Biggl({{q^{-i},\;s^*q^{i+1},\;q^{-j}}\atop
{rq,\;\;q^{-d}}}\;\Bigg\vert \; q,\;q\Biggr)
\end{eqnarray*}
for $0 \leq i,j\leq d$. These $f_i$ are $q$-Hahn polynomials.
\end{example}

\begin{example}
\label{ex:dual qhahn}
(Dual $q$-Hahn)
Assume
\begin{eqnarray*}
\theta_i &=& \theta_0+h(1-q^i)(1-sq^{i+1})q^{-i},
\\
\theta^*_i &=&\theta^*_0+h^*(1-q^i)q^{-i} 
\end{eqnarray*}
for $0 \leq i \leq d$ and
\begin{eqnarray*}
\varphi_i &=& hh^*q^{1-2i}(1-q^i)(1-q^{i-d-1})(1-rq^i),
\\
\phi_i &=& hh^*q^{d+2-2i}(1-q^i)(1-q^{i-d-1})(s-rq^{i-d-1})
\end{eqnarray*}
for $1 \leq i \leq d$.
Assume $h,h^*, q, r,s$  are nonzero.
Assume
none of
$q^i,
rq^i,
sq^i/r$
is equal to $1$ for $1 \leq i \leq d$ and 
that $sq^i\not=1$ for $2 \leq i \leq 2d$.
Then the sequence
$(\theta_i, \theta^*_i, i=0...d;  \varphi_j, \phi_j, j=1...d)$
is a parameter array over $\K$ which has base $q$.
The corresponding polynomials
$f_i$
satisfy
\begin{eqnarray*}
f_i(\theta_j)=
{}_3\phi_2 \Biggl({{q^{-i},\;q^{-j},\;sq^{j+1}}\atop
{rq,\;\;q^{-d}}}\;\Bigg\vert \; q,\;q\Biggr)
\end{eqnarray*}
for $0 \leq i,j\leq d$. These $f_i$ are dual $q$-Hahn polynomials.
\end{example}

\begin{example}
(Quantum $q$-Krawtchouk)
\label{ex:qkrawquantum}
Assume
\begin{eqnarray*}
\theta_i &=& \theta_0-sq(1-q^i),
\\
\theta^*_i &=& \theta^*_0+h^*(1-q^i)q^{-i} 
\end{eqnarray*}
for 
 $0 \leq i \leq d$  and
\begin{eqnarray*}
\varphi_i &=& -rh^*q^{1-i}(1-q^i)(1-q^{i-d-1}),
\\
\phi_i &=& h^*q^{d+2-2i}(1-q^i)(1-q^{i-d-1})(s-rq^{i-d-1})
\end{eqnarray*}
for $1 \leq i \leq d$.
Assume $ h^*, q, r, s$ are nonzero.
Assume
neither of
$q^i,
sq^i/r$
is equal to $1$ for $1 \leq i \leq d$.
Then the sequence
$(\theta_i, \theta^*_i, i=0...d;  \varphi_j, \phi_j, j=1...d)$
is a parameter array over $\K$ which has base $q$.
The corresponding polynomials
$f_i$
satisfy
\begin{eqnarray*}
f_i(\theta_j)=
{}_2\phi_1 \Biggl({{q^{-i}, \;q^{-j}}\atop
{q^{-d}}}\;\Bigg\vert \; q,\;sr^{-1}q^{j+1}\Biggr)
\end{eqnarray*}
for $0 \leq i,j\leq d$. These $f_i$ are
quantum $q$-Krawtchouk polynomials.
\end{example}

\begin{example}
($q$-Krawtchouk)
\label{ex:sbutnotss2}
\label{ex:qkraw}
Assume
\begin{eqnarray*}
\theta_i &=& 
\theta_0+h(1-q^i)q^{-i},
\\
\theta^*_i &=&  \theta^*_0+h^*(1-q^i)(1-s^*q^{i+1})q^{-i}
\end{eqnarray*}
for $0 \leq i \leq d$ and
\begin{eqnarray*}
\varphi_i &=& hh^*q^{1-2i}(1-q^i)(1-q^{i-d-1}),
\\
\phi_i &=& hh^*s^*q(1-q^i)(1-q^{i-d-1})
\end{eqnarray*}
for $1 \leq i \leq d$.
Assume $h,h^*, q, s^*$ are nonzero.
Assume  
$q^i\not=1$
for $1 \leq i \leq d$ and 
that $s^*q^i\not=1$ for $2 \leq i \leq 2d$.
Then the sequence
$(\theta_i, \theta^*_i, i=0...d;  \varphi_j, \phi_j, j=1...d)$
is a parameter array over $\K$ which has base $q$.
The corresponding polynomials
$f_i$
  satisfy
\begin{eqnarray*}
f_i(\theta_j)=
{}_3\phi_2 \Biggl({{q^{-i}, \;s^*q^{i+1},\;q^{-j}}\atop
{0,\;\;q^{-d}}}\;\Bigg\vert \; q,\;q\Biggr)
\end{eqnarray*}
for $0 \leq i,j\leq d$. These $f_i$ are $q$-Krawtchouk polynomials.
\end{example}

\begin{example}(Affine $q$-Krawtchouk)
\label{ex:qkrawaffine}
Assume
\begin{eqnarray*}
\theta_i &=& \theta_0+h(1-q^i)q^{-i},
\\
\theta^*_i &=& \theta^*_0+h^*(1-q^i)q^{-i}
\end{eqnarray*}
for 
 $0 \leq i \leq d$  and
\begin{eqnarray*}
\varphi_i &=& hh^*q^{1-2i}(1-q^i)(1-q^{i-d-1})(1-rq^i),
\\
\phi_i &=& -hh^*rq^{1-i}(1-q^i)(1-q^{i-d-1})
\end{eqnarray*}
for $1 \leq i \leq d$.
Assume $h,h^*, q, r$ are nonzero.
Assume
neither of
$q^i,
rq^i$
is equal to $1$ for $1 \leq i \leq d$.
Then the sequence
$(\theta_i, \theta^*_i, i=0...d;  \varphi_j, \phi_j, j=1...d)$
is a parameter array over $\K$ which has base $q$.
The corresponding  polynomials 
$f_i$
 satisfy
\begin{eqnarray*}
f_i(\theta_j)=
{}_3\phi_2 \Biggl({{q^{-i}, \;0,\;q^{-j}}\atop
{rq,\;\;q^{-d}}}\;\Bigg\vert \; q,\;q\Biggr)
\end{eqnarray*}
for $0 \leq i,j\leq d$. These $f_i$ are affine $q$-Krawtchouk polynomials.
\end{example}

\begin{example}
(Dual $q$-Krawtchouk)
\label{ex:sbutnotss}
Assume 
\begin{eqnarray*}
\theta_i &=& \theta_0+h(1-q^i)(1-sq^{i+1})q^{-i},
\\
\theta^*_i &=& \theta^*_0+h^*(1-q^i)q^{-i}
\end{eqnarray*}
for $0 \leq i \leq d$ and
\begin{eqnarray*}
\varphi_i &=& hh^*q^{1-2i}(1-q^i)(1-q^{i-d-1}),
\\
\phi_i &=& hh^*sq^{d+2-2i}(1-q^i)(1-q^{i-d-1})
\end{eqnarray*}
for $1 \leq i \leq d$.
Assume $h,h^*,q, s$ are nonzero.
Assume 
$q^i\not=1$
for $1 \leq i \leq d$ and 
$sq^i\not=1$ for $2 \leq i \leq 2d$.
Then the sequence
$(\theta_i, \theta^*_i, i=0...d;  \varphi_j, \phi_j, j=1...d)$
is a parameter array over $\K$ which has base $q$.
The corresponding polynomials
$f_i$
 satisfy
\begin{eqnarray*}
f_i(\theta_j)=
{}_3\phi_2 \Biggl({{q^{-i}, \;q^{-j},sq^{j+1}}\atop
{0,\;\;q^{-d}}}\;\Bigg\vert \; q,\;q\Biggr)
\end{eqnarray*}
for $0 \leq i,j\leq d$. These $f_i$ are dual $q$-Krawtchouk polynomials.
\end{example}

\begin{example}(Racah)
\label{ex:racah}
Assume 
\begin{eqnarray}
\theta_i &=& \theta_0+hi(i+1+s),
\label{eq:thracII}
\\
\theta^*_i &=& \theta^*_0+h^*i(i+1+s^*)
\label{eq:thsracII}
\end{eqnarray}
for $0 \leq i \leq d$ and
\begin{eqnarray}
\varphi_i &=& hh^*i(i-d-1)(i+r_1)(i+r_2),
\label{lem:varphiphi1II}
\\
\phi_i &=& hh^*i(i-d-1)(i+s^*-r_1)(i+s^*-r_2)
\label{lem:varphiphi2II}
\end{eqnarray}
for $1 \leq i \leq d$.
Assume $h, h^*$ are nonzero 
and that $r_1+r_2=s+s^*+d+1$.
Assume the characteristic of $\K$ is $0$ or a prime 
greater than $d$.
Assume none of
$
r_1, r_2,
s^*-r_1$, $ s^*-r_2$
is equal to $-i$ for $1 \leq i \leq d$ and 
that neither of $s, s^*$ is equal to $-i$ for $2 \leq i \leq 2d$.
Then the sequence
$(\theta_i, \theta^*_i, i=0...d;  \varphi_j, \phi_j, j=1...d)$
is a parameter array over $\K$ which has base $1$.
The corresponding polynomials
$f_i$
 satisfy
\begin{eqnarray*}
f_i(\theta_j)=
{}_4F_3 \Biggl({{-i, \;i+1+s^*,\;-j,\;j+1+s}\atop
{r_1+1,\;\;r_2+1,\;\;-d}}\;\Bigg\vert \; 1\Biggr)
\end{eqnarray*}
for $0 \leq i,j\leq d$. These $f_i$ are Racah polynomials.
\end{example}

\begin{example}(Hahn)
\label{ex:Hahn}
Assume
\begin{eqnarray*}
\theta_i &=& \theta_0+si,
\\
\theta^*_i &=& \theta^*_0+h^*i(i+1+s^*)
\end{eqnarray*}
for $0 \leq i \leq d$ and
\begin{eqnarray*}
\varphi_i &=& h^*si(i-d-1)(i+r),
\\
\phi_i &=& -h^*si(i-d-1)(i+s^*-r)
\end{eqnarray*}
for $1 \leq i \leq d$.
Assume $h^*, s$ are nonzero.
Assume the characteristic of $\K$ is $0$ or a prime greater than $d$.
Assume neither of
$
r, 
s^*-r$
is equal to $-i$ for $1 \leq i \leq d$ and 
that  $s^*\not=-i$ for $2 \leq i \leq 2d$.
Then the sequence
$(\theta_i, \theta^*_i, i=0...d;  \varphi_j, \phi_j, j=1...d)$
is a parameter array over $\K$ which has base $1$.
The corresponding polynomials
$f_i$
 satisfy
\begin{eqnarray*}
f_i(\theta_j)=
{}_3F_2 \Biggl({{-i, \;i+1+s^*,\;-j}\atop
{r+1,\;\;-d}}\;\Bigg\vert \; 1\Biggr)
\end{eqnarray*}
for $0 \leq i,j\leq d$. These $f_i$ are  Hahn polynomials.
\end{example}

\begin{example}(Dual Hahn)
\label{ex:dHahn}
Assume
\begin{eqnarray*}
\theta_i &=& \theta_0+hi(i+1+s),
\\
\theta^*_i &=& \theta^*_0+s^*i
\end{eqnarray*}
for $0 \leq i \leq d$ and
\begin{eqnarray*}
\varphi_i &=& hs^*i(i-d-1)(i+r),
\\
\phi_i &=& hs^*i(i-d-1)(i+r-s-d-1)
\end{eqnarray*}
for $1 \leq i \leq d$.
Assume $h, s^*$ are nonzero.
Assume the characteristic of $\K$ is $0$ or a prime greater than $d$.
Assume neither of
$
r, 
s-r$
is equal to $-i$ for $1 \leq i \leq d$ and 
that  $s\not=-i$ for $2 \leq i \leq 2d$.
Then the sequence
$(\theta_i, \theta^*_i, i=0...d;  \varphi_j, \phi_j, j=1...d)$
is a parameter array over $\K$ which has base $1$.
The corresponding polynomials 
$f_i$
 satisfy
\begin{eqnarray*}
f_i(\theta_j)=
{}_3F_2 \Biggl({{-i, \;-j,\;j+1+s}\atop
{r+1,\;\;-d}}\;\Bigg\vert \; 1\Biggr)
\end{eqnarray*}
for $0 \leq i,j\leq d$. These $f_i$ are dual Hahn polynomials.
\end{example}

\begin{example}(Krawtchouk)
\label{ex:kraw}
Assume
\begin{eqnarray*}
\theta_i &=& \theta_0+si,
\\
\theta^*_i &=& \theta^*_0+s^*i
\end{eqnarray*}
for $0 \leq i \leq d$ and
\begin{eqnarray*}
\varphi_i &=& ri(i-d-1)
\\
\phi_i &=& (r-ss^*)i(i-d-1)
\end{eqnarray*}
for $1 \leq i \leq d$.
Assume $r,s,s^*$ are nonzero.
Assume the characteristic of $\K$ is $0$ or a prime greater than $d$.
Assume
$
r\not=ss^*$. 
Then the sequence
$(\theta_i, \theta^*_i, i=0...d;  \varphi_j, \phi_j, j=1...d)$
is a parameter array over $\K$ which has base $1$.
The corresponding polynomials 
$f_i$
 satisfy
\begin{eqnarray*}
f_i(\theta_j)=
{}_2F_1 \Biggl({{-i, \;-j}\atop
{-d}}\;\Bigg\vert \; r^{-1}ss^*\Biggr)
\end{eqnarray*}
for $0 \leq i,j\leq d$. These $f_i$ are Krawtchouk polynomials.
\end{example}

\begin{example}(Bannai/Ito)
\label{ex:bi}
Assume
\begin{eqnarray}
\theta_i &=& \theta_0+h(s-1+(1-s+2i)(-1)^i),
\label{eq:thracIII}
\\
\theta^*_i &=& \theta^*_0+h^*(s^*-1+(1-s^*+2i)(-1)^i)
\label{eq:thsracIII}
\end{eqnarray}
for $0 \leq i \leq d$ and
\begin{eqnarray}
\varphi_i &=& 
\cases{
-4hh^*i(i+r_1),    &if $\;i$ even, $\;d$ even;\cr
-4hh^*(i-d-1)(i+r_2),    &if $\;i$ odd, $\;d$ even;\cr
-4hh^*i(i-d-1),    &if $\;i$ even, $\;d$ odd;\cr
-4hh^*(i+r_1)(i+r_2),    &if $\;i$ odd, $\;d$ odd,
}
\label{lem:varphiphi1III}
\\
\phi_i &=& 
\cases{
4hh^*i(i-s^*-r_1),    &if $\;i$ even, $\;d$ even;\cr
4hh^*(i-d-1)(i-s^*-r_2),    &if $\;i$ odd, $\;d$ even;\cr
-4hh^*i(i-d-1),    &if $\;i$ even, $\;d$ odd;\cr
-4hh^*(i-s^*-r_1)(i-s^*-r_2),    &if $\;i$ odd, $\;d$ odd
}
\label{lem:varphiphi2III}
\end{eqnarray}
for $1 \leq i \leq d$.
Assume $h, h^*$ are nonzero and that
$r_1+r_2=-s-s^*+d+1$.
Assume the characteristic of $\K$ is either $0$
or an odd prime greater
than $d/2$.
Assume neither of
$
r_1, -s^*-r_1$
is equal to $-i$ for $1 \leq i \leq d$, $d-i$ even.
Assume neither of 
$r_2$, $ -s^*-r_2$
is equal to $-i$ for $1 \leq i \leq d$, $i$ odd.
Assume neither of $s, s^*$ is equal to $2i$ for $1 \leq i \leq d$.
Then the sequence
$(\theta_i, \theta^*_i, i=0...d;  \varphi_j, \phi_j, j=1...d)$
is a parameter array over $\K$ which has base $-1$.
We call the corresponding polynomials  from
Definition \ref{def:assocpoly}
the Bannai/Ito polynomials 
\cite[p. 260]{BanIto}.
\end{example}

\begin{example} (Orphan) 
\label{ex:orphan}
For this example assume
$\K$ has characteristic 2.
For notational convenience we
define some scalars $\gamma_0, 
\gamma_1, 
\gamma_2, 
\gamma_3$  
in $ \K$.
We define
$\gamma_i=0$ for $i \in \lbrace 0,3\rbrace $
 and $\gamma_i=1$ for $i\in \lbrace 1,2\rbrace$.
Assume
\begin{eqnarray}
\theta_i &=& \theta_0 + h(si+ \gamma_i),
\label{eq:thracIV}
\\
\theta^*_i &=& \theta^*_0 + h^*(s^*i+ \gamma_i)
\label{eq:thsracIV}
\end{eqnarray}
for $0 \leq i \leq 3$.
Assume $\varphi_1 = hh^*r$, 
$\varphi_2 = hh^*$, 
$\varphi_3 = hh^*(r+s+s^*)$ 
and 
$\phi_1 = hh^*(r+s(1+s^*))$, 
$\phi_2 = hh^*$, 
$\phi_3 = hh^*(r+s^*(1+s))$.
Assume each of $h,h^*, s,s^*,r$ is nonzero.
Assume neither of $s, s^*$ is equal to $1$
and that $r$ is equal to none of 
$s+s^*$, 
$s(1+s^*)$, 
$s^*(1+s)$.
Then the sequence
$(\theta_i, \theta^*_i, i=0...3;  \varphi_j, \phi_j, j=1...3)$
is a parameter array over $\K$ which has
diameter $3$ and base $1$.
We call the corresponding polynomials  from 
Definition \ref{def:assocpoly}
the Orphan polynomials.
\end{example}


\begin{theorem}
\label{thm:thatsit}
Every parameter array over $\K$
is listed in at least one of the Examples
\ref{ex:qracah}--\ref{ex:orphan}.
\end{theorem}
\noindent {\it Proof:}
Let 
$p:=(\theta_i, \theta^*_i, i=0...d;  \varphi_j, \phi_j, j=1...d)$
denote a parameter array over $\K$.
We show this array is given in
at least one of the Examples \ref{ex:qracah}--\ref{ex:orphan}.
We assume $d\geq 1$; otherwise the result is trivial. 
For notational convenience let ${\tilde K}$ denote the 
algebraic closure of $\K$.
Let $q$ denote a base for $p$ as in Definition
\ref{def:base}.
For $d<3$ we may assume
$q\not=1$ and $q \not=-1$ in view of the remark in 
Definition
\ref{def:base}.
By PA5 and 
Definition 
\ref{def:base}
both
\begin{eqnarray}
\theta_{i-2}- \lbrack 3 \rbrack_q\theta_{i-1}+
\lbrack 3 \rbrack_q\theta_{i} -\theta_{i+1}&=&0,
\label{eq:lin1}
\\
\theta^*_{i-2}- \lbrack 3 \rbrack_q\theta^*_{i-1}+
\lbrack 3 \rbrack_q\theta^*_{i} -\theta^*_{i+1}&=&0
\label{eq:lin2}
\end{eqnarray}
for $2 \leq i \leq d-1$, where $\lbrack 3 \rbrack_q := q+q^{-1}+1$.
We divide the argument into the following four cases.
(I) $q \not=1$, $q \not=-1$;
(II) $q = 1$ and $\mbox{char}(\K) \not=2$;
(III) $q = -1$ and $\mbox{char}(\K) \not=2$;
(IV) $q= 1 $ and 
 $\mbox{char}(\K)=2$.

\medskip
\noindent Case I: $q\not=1$, $q\not=-1$.  \\
\noindent By (\ref{eq:lin1})  there exist scalars $\eta, \mu, h$ in
${\tilde \K}$ such that
\begin{eqnarray}
\theta_i &=& \eta + \mu q^i+ h q^{-i}  \qquad \qquad 
(0 \leq i\leq d).
\label{eq:thI}
\end{eqnarray}
By 
(\ref{eq:lin2}) there exist 
 scalars $\eta^*, \mu^*,h^*$ in
${\tilde \K}$ such that
\begin{eqnarray}
\theta^*_i &=& \eta^* + \mu^* q^i+ h^*q^{-i}  \qquad \qquad 
(0 \leq i\leq d).
\label{eq:thsI}
\end{eqnarray}
Observe $\mu, h$ are not both 0; 
otherwise $\theta_1=\theta_0$ by
(\ref{eq:thI}).
Similarly 
$\mu^*, h^*$ are not both 0.
For $1 \leq i \leq d$ we have $q^i \not=1$; otherwise
$\theta_i=\theta_0$ by
(\ref{eq:thI}).
Setting $i=0$ 
in 
(\ref{eq:thI}), (\ref{eq:thsI}) we obtain
\begin{eqnarray}
\theta_0 &=& \eta + \mu +h,
\label{eq:thI0}
\\
\theta^*_0 &=& \eta^* + \mu^*+h^*.
\label{eq:thsI0}
\end{eqnarray}
We claim there exists $\tau \in {\tilde \K}$ such that both
\begin{eqnarray}
\varphi_i &=&
(q^i-1)(q^{d-i+1}-1)(\tau -\mu \mu^* q^{i-1} - h h^* q^{-i-d}),
\label{eq:varphiIfirst}
\\
\phi_i &=&
(q^i-1)(q^{d-i+1}-1)(\tau -h \mu^* q^{i-d-1} - \mu h^* q^{-i})
\label{eq:phiIfirst}
\end{eqnarray}
for $1 \leq i \leq d$.
Since $q\not=1$ and $q^d\not=1$ there exists $\tau \in {\tilde \K}$
such that 
(\ref{eq:varphiIfirst})
holds for $i=1$.
In the equation of PA4, we eliminate
$\varphi_1$  using  
(\ref{eq:varphiIfirst}) at $i=1$, and evaluate the result using
(\ref{eq:thI}), 
(\ref{eq:thsI}), and
\cite[Lemma 10.2]{LS99} in order to obtain 
(\ref{eq:phiIfirst}) for $1 \leq i \leq d$.
In the equation of PA3, we
eliminate $\phi_1$ using  
(\ref{eq:phiIfirst}) at $i=1$, and evaluate the result using
(\ref{eq:thI}), 
(\ref{eq:thsI}),
and \cite[Lemma 10.2]{LS99} in order to obtain
(\ref{eq:varphiIfirst})  for $1 \leq i \leq d$.
We have now proved the claim. 
We now break the argument into subcases.
For each subcase our argument is similar.
We will discuss the first subcase in detail in order to
give the idea; for the remaining subcases 
we give the essentials only.

\medskip
\noindent Subcase $q$-Racah: 
$\mu\not=0, \mu^*\not=0, h\not=0, h^*\not=0$.  
We show $p$ is listed in Example
\ref{ex:qracah}.
Define 
\begin{eqnarray}
s:=\mu h^{-1} q^{-1}, \qquad   
s^*:=\mu^* h^{*-1} q^{-1}.
\label{eq:hdef}
\end{eqnarray}
 Eliminating
$\eta $ in 
(\ref{eq:thI})  using
(\ref{eq:thI0}) and eliminating $\mu $ in the result
using the equation on the left in  
(\ref{eq:hdef}),
we obtain
(\ref{eq:thrac}) for $0 \leq i\leq d$.
Similarly we obtain
(\ref{eq:thsrac}) for $0 \leq i\leq d$.
Since $\tilde \K $ is algebraically closed it contains
scalars $r_1, r_2$
such that both
\begin{eqnarray}
r_1r_2 = s s^*q^{d+1}, \qquad \qquad 
r_1+r_2=\tau h^{-1} h^{*-1} q^d.
\label{eq:rdef}
\end{eqnarray}
Eliminating $\mu, \mu^*, \tau$ in
(\ref{eq:varphiIfirst}), 
(\ref{eq:phiIfirst}) using
(\ref{eq:hdef}) and the equation on the right in
(\ref{eq:rdef}), and evaluating the result using
the equation on the left in 
(\ref{eq:rdef}), 
we obtain 
(\ref{lem:varphiphi1}), 
(\ref{lem:varphiphi2}) 
for $1 \leq i \leq d$.
By the construction 
each of $h, h^*, q, s, s^*$ is nonzero.
Each of $r_1, r_2$ is nonzero by the equation on the left in
(\ref{eq:rdef}). The remaining inequalities mentioned below
(\ref{lem:varphiphi2}) follow
from PA1, PA2 and 
(\ref{eq:thrac})--(\ref{lem:varphiphi2}).
We have now shown $p$ is listed in Example
\ref{ex:qracah}.

\medskip
\noindent
We now give the remaining subcases of Case I.
We list the essentials only. 

\medskip
\noindent Subcase $q$-Hahn:
$\mu=0, \mu^*\not=0, h\not=0, h^*\not=0, \tau \not=0$.  
Definitions:
\begin{eqnarray*}
s^*:=\mu^* h^{*-1} q^{-1},\qquad
r:=\tau h^{-1} h^{*-1} q^d.
\end{eqnarray*}
\noindent Subcase dual $q$-Hahn:
$\mu\not=0, \mu^*=0, h\not=0, h^*\not=0, \tau \not=0$.  
Definitions:
\begin{eqnarray*}
 s:=\mu h^{-1} q^{-1},\qquad
r:=\tau h^{-1} h^{*-1} q^d.
\end{eqnarray*}
\noindent Subcase quantum $q$-Krawtchouk: 
$\mu\not=0, \mu^*=0, h=0, h^*\not=0, \tau \not=0$. 
Definitions: 
\begin{eqnarray*}
s:=\mu  q^{-1}, \qquad   
r :=\tau h^{*-1}q^d.
\end{eqnarray*}
\noindent Subcase  $q$-Krawtchouk: 
$\mu=0, \mu^*\not=0, h\not=0, h^*\not=0, \tau =0$. 
Definition: 
\begin{eqnarray*}
s^*:=
 \mu^* h^{*-1} q^{-1}.
\end{eqnarray*}
\noindent Subcase affine $q$-Krawtchouk:
$\mu=0, \mu^*=0, h\not=0, h^*\not=0, \tau \not=0$. 
Definition:
\begin{eqnarray*}
r:=\tau h^{-1} h^{*-1} q^d.
\end{eqnarray*}
\noindent Subcase dual $q$-Krawtchouk:
$\mu \not=0, \mu^*=0, h\not=0, h^*\not=0, \tau =0$. 
Definition:
\begin{eqnarray*}
s:=
 \mu h^{-1} q^{-1}.
\end{eqnarray*}
\noindent 
We have a few more comments concerning Case I. 
Earlier we mentioned that $\mu,h$ are not both 0 and 
that $\mu^*,h^*$ are not both 0.
Suppose one of 
$\mu,h$ is 0 
and one of  $\mu^*,h^*$ is 0.
Then $\tau\not=0$; otherwise $\varphi_1=0$ by 
(\ref{eq:varphiIfirst}) or 
$\phi_1=0$ by 
(\ref{eq:phiIfirst}).
Suppose  $\mu^*\not=0$,  $h^*=0$.
Replacing $q$ by $q^{-1}$
we obtain  
$\mu^* =0$, $h^*\not=0$. 
Suppose 
 $\mu^* \not=0$, $h^*\not=0$,  
 $\mu\not=0$,  $h=0$. Replacing $q$ by $q^{-1}$
we obtain
 $\mu^* \not=0$, $h^*\not=0$,  
 $\mu=0$,  $h\not=0$.
By these comments we find that after 
replacing $q$ by $q^{-1}$ if necessary,
one of the above subcases holds.
This completes our argument for Case I.

\medskip
\noindent Case II: $q=1$ and $\mbox{char}(\K)\not=2$.
\\
\noindent By (\ref{eq:lin1}) and since
$\mbox{char}(\K)\not=2$, 
there exist scalars $\eta, \mu, h$ in
${\tilde \K}$ such that
\begin{eqnarray}
\theta_i &=& \eta + (\mu + h) i+ h i^2  \qquad \qquad 
(0 \leq i\leq d).
\label{eq:thII}
\end{eqnarray}
Similarly  
there exist 
 scalars $\eta^*, \mu^*,h^*$ in
${\tilde \K}$ such that
\begin{eqnarray}
\theta^*_i &=& \eta^* + (\mu^* +h^*)i+ h^*i^2  \qquad \qquad 
(0 \leq i\leq d).
\label{eq:thsII}
\end{eqnarray}
Observe $\mu, h$ are not both 0; 
otherwise $\theta_1=\theta_0$.
Similarly 
$\mu^*, h^*$ are not both 0.
For any prime $i$ such that $ i \leq d$ we have
$\mbox{char}(\K)\not=i$;
otherwise
$\theta_i=\theta_0$ by
(\ref{eq:thII}).
Therefore
$\mbox{char}(\K)$ is 0 or a prime greater than $d$.
Setting $i=0$ 
in 
(\ref{eq:thII}), (\ref{eq:thsII}) we obtain
\begin{eqnarray}
\theta_0 = \eta, \qquad \qquad 
\theta^*_0 = \eta^*.
\label{eq:thII0}
\end{eqnarray}
We claim there exists $\tau \in {\tilde \K}$ such that both
\begin{eqnarray}
\varphi_i &=&
i(d-i+1)(\tau -( \mu h^*+h \mu^*) i - h h^* i(i+d+1)),
\label{eq:varphiIIfirst}
\\
\phi_i &=&
i(d-i+1)(\tau +\mu \mu^*+ h \mu^*(1+d)  +(\mu h^*-h \mu^*) i + h h^* i(d-i+1))
\label{eq:phiIIfirst}
\end{eqnarray}
for $1 \leq i \leq d$.
There exists $\tau \in {\tilde \K}$
such that 
(\ref{eq:varphiIIfirst})
holds for $i=1$.
In the equation of PA4, we eliminate
$\varphi_1$  using  
(\ref{eq:varphiIIfirst}) at $i=1$, and evaluate the result using
(\ref{eq:thII}), 
(\ref{eq:thsII}), and
\cite[Lemma 10.2]{LS99} in order to obtain 
(\ref{eq:phiIIfirst}) for $1 \leq i \leq d$.
In the equation of PA3, we
eliminate $\phi_1$ using  
(\ref{eq:phiIIfirst}) at $i=1$, and evaluate the result using
(\ref{eq:thII}), 
(\ref{eq:thsII}),
and \cite[Lemma 10.2]{LS99} in order to obtain
(\ref{eq:varphiIIfirst})  for $1 \leq i \leq d$.
We have now proved the claim. 
We now break the argument into subcases.

\medskip
\noindent Subcase Racah: 
$h\not=0$, $h^*\not=0$.  
We show $p$ is listed in Example
\ref{ex:racah}.
Define 
\begin{eqnarray}
s:=\mu h^{-1} , \qquad   
s^*:=\mu^* h^{*-1}.
\label{eq:hdefII}
\end{eqnarray}
 Eliminating
$\eta, \mu $ in 
(\ref{eq:thII})  using
(\ref{eq:thII0}),
(\ref{eq:hdefII})
we obtain
(\ref{eq:thracII}) for $0 \leq i\leq d$.
Eliminating
$\eta^*, \mu^*$  in 
(\ref{eq:thsII})  using
(\ref{eq:thII0}),
(\ref{eq:hdefII})
we obtain
(\ref{eq:thsracII}) for $0 \leq i\leq d$.
Since $\tilde \K $ is algebraically closed it contains
scalars $r_1, r_2$
such that both
\begin{eqnarray}
r_1r_2 =  - \tau h^{-1} h^{*-1}, \qquad \qquad 
r_1+r_2=s + s^*+d+1.
\label{eq:rdefII}
\end{eqnarray}
Eliminating $\mu, \mu^*, \tau$ in
(\ref{eq:varphiIIfirst}), 
(\ref{eq:phiIIfirst}) using
(\ref{eq:hdefII}) and the equation on the left in
(\ref{eq:rdefII})
we obtain 
(\ref{lem:varphiphi1II}), 
(\ref{lem:varphiphi2II}) 
for $1 \leq i \leq d$.
By the construction 
each of $h, h^*$ is nonzero.
The remaining inequalities mentioned below
(\ref{lem:varphiphi2II}) follow
from PA1, PA2 and 
(\ref{eq:thracII})--(\ref{lem:varphiphi2II}).
We have now shown $p$ is listed in Example
\ref{ex:racah}.

\medskip
\noindent We now give the remaining subcases of Case II. We list the
essentials only.

\medskip
\noindent Subcase Hahn:
$ h=0$, $h^*\not=0 $.  
Definitions:
\begin{eqnarray*}
s= \mu, 
\qquad s^*:=\mu^* h^{*-1}, \qquad 
r:=-\tau \mu^{-1} h^{*-1}.
\end{eqnarray*}
\noindent Subcase dual Hahn:
$ h\not=0, h^*=0$.  
Definitions:
\begin{eqnarray*}
 s:=\mu h^{-1},\qquad
s^*= \mu^*, \qquad 
r:=- \tau 
h^{-1}
\mu^{*-1} .
\end{eqnarray*}
\noindent Subcase Krawtchouk: 
$ h=0, h^*=0$. 
Definitions: 
\begin{eqnarray*}
s:=\mu , \qquad   
s^*:=\mu^* , \qquad   
r :=- \tau.
\end{eqnarray*}

\medskip
\noindent Case III: $q = -1$ and $\mbox{char}(\K)\not=2$.
\\
\noindent
We show $p$ is listed in
Example \ref{ex:bi}.
By (\ref{eq:lin1})  and since
$\mbox{char}(\K)\not=2$, 
there exist scalars $\eta, \mu, h$ in
${\tilde \K}$ such that
\begin{eqnarray}
\theta_i &=& \eta + \mu (-1)^i+ 2h i(-1)^i  \qquad \qquad 
(0 \leq i\leq d).
\label{eq:thIII}
\end{eqnarray}
Similarly  
there exist 
 scalars $\eta^*, \mu^*,h^*$ in
${\tilde \K}$ such that
\begin{eqnarray}
\theta^*_i &=& \eta^* + \mu^*(-1)^i+ 2 h^*i(-1)^i  \qquad \qquad 
(0 \leq i\leq d).
\label{eq:thsIII}
\end{eqnarray}
Observe $h\not=0$;
otherwise $\theta_2=\theta_0$ by
(\ref{eq:thIII}).
Similarly 
$h^*\not=0$.
For any prime $i$ such that $ i \leq d/2$ we have
$\mbox{char}(\K)\not=i$;
otherwise
$\theta_{2i}=\theta_0$ by
(\ref{eq:thIII}).
By this and since 
$\mbox{char}(\K)\not=2$
we find $\mbox{char}(\K)$ is either
0 or an odd prime greater than $d/2$.
Setting $i=0$ 
in 
(\ref{eq:thIII}), (\ref{eq:thsIII}) we obtain
\begin{eqnarray}
\theta_0 &=& \eta+\mu, 
\qquad \qquad 
\theta^*_0 = \eta^*+\mu^*.
\label{eq:thIII0}
\end{eqnarray}
We define 
\begin{eqnarray}
s:=1-\mu h^{-1}, \qquad \qquad  s^*=1-\mu^*h^{*-1}. 
\label{eq:sssIII}
\end{eqnarray}
Eliminating $\eta$ in 
(\ref{eq:thIII}) using
(\ref{eq:thIII0}) and eliminating $\mu$ in the result using
(\ref{eq:sssIII}) we find 
(\ref{eq:thracIII}) holds for $0 \leq i \leq d$.
Similarly we find
(\ref{eq:thsracIII}) holds for $0 \leq i \leq d$.
We now define $r_1, r_2$. First 
assume 
$d$ is odd. 
Since $\tilde \K$ is algebraically closed it contains 
$r_1, r_2$ such that
\begin{eqnarray}
r_1+r_2=-s-s^*+d+1
\label{r1pr2III}
\end{eqnarray}
and such that
\begin{eqnarray}
4hh^*(1+r_1)(1+r_2)= -\varphi_1.
\label{eq:rprodIII}
\end{eqnarray}
Next assume $d$ is even. Define
\begin{eqnarray}
r_2 := -1 + \frac{\varphi_1}{4hh^*d}
\label{eq:deveIII}
\end{eqnarray}
and define $r_1$ so that 
(\ref{r1pr2III}) holds.
We have now defined $r_1, r_2$ for either parity of $d$.
In the equation of PA4, we eliminate
$\varphi_1$  using 
(\ref{eq:rprodIII}) or 
(\ref{eq:deveIII}), and evaluate the result using
(\ref{eq:thracIII}), 
(\ref{eq:thsracIII}), 
and \cite[Lemma 10.2]{LS99} in order to obtain
(\ref{lem:varphiphi2III}) for $1 \leq i\leq d$.
In the equation of PA3, we
eliminate $\phi_1$ using  
(\ref{lem:varphiphi2III}) 
at $i=1$, and evaluate the result using
(\ref{eq:thracIII}), 
(\ref{eq:thsracIII}),
and \cite[Lemma 10.2]{LS99} in order to obtain
(\ref{lem:varphiphi1III}) for $1 \leq i\leq d$.
We mentioned each of $h, h^*$ is nonzero.
The remaining inequalities mentioned below
(\ref{lem:varphiphi2III}) 
follow from PA1, PA2 and
(\ref{eq:thracIII})--(\ref{lem:varphiphi2III}).
We have now shown $p$  is listed in
Example \ref{ex:bi}.

\medskip
\noindent Case IV: $q = 1$ and $\mbox{char}(\K)=2$.
\\
\noindent 
We show $p$ is listed in
  Example
\ref{ex:orphan}.
We first show $d=3$.
Recall $d\geq 3$
since $q=1$.
Suppose $d\geq 4$.
By 
(\ref{eq:lin1}) we have $\sum_{j=0}^3 \theta_j=0$
and 
$\sum_{j=1}^4 \theta_j=0$.
Adding 
these sums
we find $\theta_0=\theta_4$ which contradicts
PA1.
Therefore $d=3$. 
We claim there exist nonzero scalars 
$h,s$ in $\K$ such that 
(\ref{eq:thracIV}) holds for $0 \leq i \leq 3$.
Define $h=\theta_0+\theta_2$. Observe $h\not=0$; otherwise
$\theta_0=\theta_2$.
Define $s=(\theta_0+\theta_3)h^{-1}$.
Observe
$s \not=0$; otherwise $\theta_0=\theta_3$.
Using these values for $h,s$ we find 
(\ref{eq:thracIV}) holds for $i=0,2,3$.
By this and $\sum_{j=0}^3 \theta_j=0$ we find 
(\ref{eq:thracIV}) holds for $i=1$.
We have now proved our claim.
Similarly there exist nonzero scalars
$h^*,s^*$ in $\K$ such that 
(\ref{eq:thsracIV}) holds for $0 \leq i \leq 3$.
Define 
$r:=\varphi_1h^{-1}h^{*-1}$.
Observe $r\not=0$ and that
$\varphi_1=h h^*r$.
In the equation of PA4, we eliminate
$\varphi_1$  using  
$\varphi_1=h h^*r$
 and evaluate the result using
(\ref{eq:thracIV}),
(\ref{eq:thsracIV})
 and
\cite[Lemma 10.2]{LS99} in order to obtain 
$\phi_1 = hh^*(r+s(1+s^*))$, 
$\phi_2 = hh^*$, 
$\phi_3 = hh^*(r+s^*(1+s))$.
In the equation of PA3, we
eliminate $\phi_1$ using  
$\phi_1 = hh^*(r+s(1+s^*))$ 
 and evaluate the result using
(\ref{eq:thracIV}),
(\ref{eq:thsracIV})
and \cite[Lemma 10.2]{LS99} in order to obtain
$\varphi_2 = hh^*$, 
$\varphi_3 = hh^*(r+s+s^*)$.
We mentioned each of $h, h^*, s,s^*, r$ is nonzero.
Observe $s\not=1$; otherwise $\theta_1=\theta_0$.
Similarly 
 $s^*\not=1$.
Observe $r\not=s+s^*$; otherwise $\varphi_3=0$.
Observe $r\not=s(1+s^*)$; otherwise $\phi_1=0$.
Observe $r\not=s^*(1+s)$; otherwise $\phi_3=0$.
We have now shown  $p$ is
listed in Example
\ref{ex:orphan}. 
We are done with Case IV and the proof is complete.
\hfill $\Box $ \\

\section{The orthogonality relation in terms of the parameter array}

Some facts about the polynomials
in Examples
\ref{ex:qracah}--\ref{ex:orphan}
can be expressed 
in a uniform and attractive manner
by writing things in terms of 
the associated parameter array.
We illustrate this by  giving the
orthogonality relation, the three-term recurrence,
and the difference equation in terms of the parameter array.
We start with the orthogonality relation. 
In order to state the result
we define some scalars.

\begin{definition}
\label{def:ki}
Let 
$(\theta_i, \theta^*_i,i=0...d; \varphi_j, \phi_j,j=1...d)$
denote a parameter array over $\K$. 
For $0 \leq i \leq d$ we let $k_i$ equal
\begin{eqnarray*}
\frac{\varphi_1\varphi_2\cdots \varphi_i}
{\phi_1\phi_2\cdots \phi_i}
\end{eqnarray*}
\noindent times 
\begin{eqnarray*}
\frac{
(\theta^*_0-\theta^*_1)
(\theta^*_0-\theta^*_2) \cdots 
(\theta^*_0-\theta^*_d)
}
{
(\theta^*_i-\theta^*_0) \cdots 
(\theta^*_i-\theta^*_{i-1}) 
(\theta^*_i-\theta^*_{i+1}) \cdots 
(\theta^*_i-\theta^*_d)
}.
\end{eqnarray*}
For $0 \leq i \leq d$ we let $k^*_i$ equal
\begin{eqnarray*}
\frac{\varphi_1\varphi_2\cdots \varphi_i}
{\phi_d\phi_{d-1}\cdots \phi_{d-i+1}}
\end{eqnarray*}
\noindent times 
\begin{eqnarray*}
\frac{
(\theta_0-\theta_1)
(\theta_0-\theta_2) \cdots 
(\theta_0-\theta_d)
}
{
(\theta_i-\theta_0) \cdots 
(\theta_i-\theta_{i-1}) 
(\theta_i-\theta_{i+1}) \cdots 
(\theta_i-\theta_d)
}.
\end{eqnarray*}
\noindent We observe $k_0=1$, $k^*_0=1$.  We define
\begin{eqnarray*}
\nu = \frac
{
(\theta_0-\theta_1)
(\theta_0-\theta_2)
\cdots
(\theta_0-\theta_d)
(\theta^*_0-\theta^*_1)
(\theta^*_0-\theta^*_2)
\cdots 
(\theta^*_0-\theta^*_d)
}
{
\phi_1\phi_2 \cdots \phi_d
}.
\end{eqnarray*}
\end{definition}

\begin{theorem} 
 \cite[Lines (128), (129)]{LS24}.
\label{thm:orth}
Let 
$(\theta_i, \theta^*_i,i=0...d; \varphi_j, \phi_j,j=1...d)$
denote a parameter array over $\K$ and 
let $f_0, f_1, \ldots, f_d$ 
  denote the corresponding polynomials
from 
Definition
\ref{def:assocpoly}.
Then both
\begin{eqnarray}
\sum_{r=0}^d f_i(\theta_r)f_j(\theta_r)k^*_r &=& \delta_{ij}\nu k^{-1}_i
		\qquad \qquad  (0 \leq i,j\leq d),
\label{eq:fiorthog}
\\
\sum_{r=0}^d f_r(\theta_i)f_r(\theta_j)k_r &=& \delta_{ij}\nu k^{*-1}_i
		\qquad \qquad  (0 \leq i,j\leq d).
\label{eq:fisorthog}
\end{eqnarray}
The scalars $k_i, k^*_i, \nu$ are from
Definition
\ref{def:ki}.
\end{theorem}

\noindent We have a comment.

\begin{lemma}
\label{lem:ksum}
With reference to Definition
\ref{def:ki}, both 
\begin{eqnarray}
\nu = \sum_{r=0}^d k_r,
\qquad \qquad \qquad 
\nu = \sum_{r=0}^d k^*_r.
\label{eq:ksum}
\end{eqnarray}
\end{lemma}
\noindent {\it Proof:}
To get the equation on the left in 
(\ref{eq:ksum}) set $i=0, j=0$ in
(\ref{eq:fisorthog}) and observe $f_r(\theta_0)=1$ for $0 \leq r \leq d$.
To get the equation on the right in  
(\ref{eq:ksum}) set $i=0, j=0$ in
(\ref{eq:fiorthog}) and observe $f_0=1$.
\hfill $\Box $ \\

\section{The three-term recurrence in terms of the parameter array}

\noindent In this section we give a three-term recurrence
satisfied by the polynomials
in Example
\ref{ex:qracah}--\ref{ex:orphan}. We express the result in
terms of the associated parameter array.
In order to state the result we define some scalars.

\begin{definition}
\label{def:bici}
Let 
$(\theta_i, \theta^*_i, i=0...d;  \varphi_j, \phi_j, j=1...d)$
denote a parameter array over $\K$.
We define
\begin{eqnarray}
b_i = \varphi_{i+1}\frac{(\theta^*_i-\theta^*_0)
(\theta^*_i-\theta^*_1) \cdots 
(\theta^*_i-\theta^*_{i-1})}
{(\theta^*_{i+1}-\theta^*_0)
(\theta^*_{i+1}-\theta^*_1) \cdots 
(\theta^*_{i+1}-\theta^*_i)} 
\qquad \qquad (0 \leq i \leq d-1) 
\label{biorig}
\end{eqnarray}
and $b_d=0$. We define
\begin{eqnarray}
c_i = \phi_{i}\frac{(\theta^*_i-\theta^*_d)
(\theta^*_i-\theta^*_{d-1}) \cdots 
(\theta^*_i-\theta^*_{i+1})}
{(\theta^*_{i-1}-\theta^*_d)
(\theta^*_{i-1}-\theta^*_{d-1}) \cdots 
(\theta^*_{i-1}-\theta^*_i)}
\qquad \qquad (1 \leq i \leq d) 
\label{ciorig}
\end{eqnarray}
and $c_0=0$. We define
\begin{eqnarray}
a_i = \theta_0-c_i - b_i \qquad \qquad (0 \leq i \leq d).
\label{aiorig}
\end{eqnarray}
\end{definition}

\begin{theorem}
\label{thm:rec1}
Let 
$(\theta_i, \theta^*_i, i=0...d;  \varphi_j, \phi_j, j=1...d)$
denote a parameter array over $\K$ and 
let  $f_0, f_1, \ldots, f_d$
denote the corresponding
polynomials from
Definition
\ref{def:assocpoly}.
For $0 \leq i,j\leq d$ 
we have
\begin{eqnarray}
\theta_j f_i(\theta_j) = c_if_{i-1}(\theta_j)
+
a_i f_i(\theta_j)
+ b_i
f_{i+1}(\theta_j),
\label{eq:threeterm}
\end{eqnarray}
where $f_{-1}, f_{d+1}$ are indeterminates
and where the $a_i, b_i, c_i$ are
from 
Definition
\ref{def:bici}.
\end{theorem}
\noindent {\it Proof:}
Let $\mathcal P$ denote the matrix in 
$\hbox{Mat}_{d+1}(\K)$
which has $ij$th entry $f_j(\theta_i)$ for $0 \leq i,j \leq  d$.
Let $K^*$ denote the diagonal matrix in
$\hbox{Mat}_{d+1}(\K)$ which has entries 
$K^*_{ii}= k^*_i$ for $0 \leq i \leq d$, where the $k^*_i$
are from
Definition
\ref{def:ki}.
Let $H$ denote the diagonal matrix in
$\hbox{Mat}_{d+1}(\K)$ which has entries 
$H_{ii}= \theta_i$ for $0 \leq i \leq d$.
Let $C$ denote the following matrix in 
$\hbox{Mat}_{d+1}(\K)$.
$$
C=\left(
\begin{array}{c c c c c c}
a_0 & b_0& & & & {\mathbf 0} \\
c_1 & a_1 & b_1  & & & \\
& c_2 & \cdot & \cdot & & \\
& & \cdot & \cdot & \cdot  &  \\
& & & \cdot & \cdot & b_{d-1} \\
{\mathbf 0}& & & & c_d & a_d 
\end{array}
\right).
$$
We define $P^*={\mathcal P}^t K^*$.
By \cite[Line (118)]{LS24} we have
$CP^*=P^*H$.
By this and since $K^*, H$ are diagonal we find
$C {\mathcal P}^t = {\mathcal P}^t H$.
In this equation we expand each side using 
matrix multiplication and routinely obtain
(\ref{eq:threeterm}).
\hfill $\Box $ \\

\noindent 
We finish this section with a comment.

\begin{lemma}
\label{lem:kialt}
With reference to Definition
\ref{def:ki} and 
Definition
\ref{def:bici}, 
\begin{eqnarray*}
k_i = \frac{
b_0b_1\cdots b_{i-1}
}
{
c_1c_2\cdots c_i
}
\qquad \qquad (0 \leq i \leq d).
\end{eqnarray*}
\end{lemma}
\noindent {\it Proof:}
Compare the formulae for the $k_i, b_i, c_i$ given
in Definition
\ref{def:ki} and 
Definition
\ref{def:bici}.
\hfill $\Box $ \\

\section{The difference equation in terms of the parameter array}

\noindent In this section we give a 
difference equation satisfied by the polynomials
in Example
\ref{ex:qracah}--\ref{ex:orphan}. We express the result in
terms of the associated parameter array.
In order to state the result we define some scalars.

\begin{definition}
\label{def:bicidual}
Let 
$(\theta_i, \theta^*_i, i=0...d;  \varphi_j, \phi_j, j=1...d)$
denote a parameter array over $\K$.
We define
\begin{eqnarray*}
b^*_i = \varphi_{i+1}\frac{(\theta_i-\theta_0)
(\theta_i-\theta_1) \cdots 
(\theta_i-\theta_{i-1})}
{(\theta_{i+1}-\theta_0)
(\theta_{i+1}-\theta_1) \cdots 
(\theta_{i+1}-\theta_i)}
\qquad \qquad (0\leq i \leq d-1)
\end{eqnarray*}
and $b^*_d=0$. We define
\begin{eqnarray*}
c^*_i = \phi_{d-i+1}\frac{(\theta_i-\theta_d)
(\theta_i-\theta_{d-1}) \cdots 
(\theta_i-\theta_{i+1})}
{(\theta_{i-1}-\theta_d)
(\theta_{i-1}-\theta_{d-1}) \cdots 
(\theta_{i-1}-\theta_i)}
\qquad \qquad (1\leq i \leq d)
\end{eqnarray*}
and $c^*_0=0$.
We define
\begin{eqnarray*}
a^*_i = \theta^*_0-c^*_i-b^*_i \qquad \qquad (0 \leq i\leq d).
\end{eqnarray*}
\end{definition}

\begin{theorem}
Let 
$(\theta_i, \theta^*_i, i=0...d;  \varphi_j, \phi_j, j=1...d)$
denote a parameter array over $\K$ and let
  $f_0, f_1, \ldots, f_d$
denote the corresponding
polynomials from
Definition 
\ref{def:assocpoly}.
For $0 \leq i,j\leq d$ we have
\begin{eqnarray}
\theta^*_if_i(\theta_j) =
c^*_jf_{i}(\theta_{j-1})
+
a^*_jf_i(\theta_j)
+ b^*_j
f_{i}(\theta_{j+1}),
\label{eq:dualrec}
\end{eqnarray}
where $\theta_{-1}, \theta_{d+1}$ are indeterminates
and
the $a^*_j, b^*_j, c^*_j $ are from
Definition \ref{def:bicidual}.
\end{theorem}
\noindent {\it Proof:}
By Lemma \ref{lem:d4}(i)
the sequence
$
(\theta^*_i, \theta_i, i=0...d;  \varphi_j, \phi_{d-j+1}, j=1...d)$
is a parameter array over $\K$.
Let  $f^*_0, f^*_1, \ldots, f^*_d$ denote
the corresponding polynomials from 
Definition
\ref{def:assocpoly},
so that
\begin{eqnarray}
f^*_i=\sum_{n=0}^i \frac{
(\lambda-\theta^*_0)
(\lambda-\theta^*_1) \cdots
(\lambda-\theta^*_{n-1})
(\theta_i-\theta_0)
(\theta_i-\theta_1) \cdots
(\theta_i-\theta_{n-1})
}
{\varphi_1\varphi_2\cdots \varphi_n}
\label{eq:dualpoly1}
\end{eqnarray}
for $0 \leq i \leq d$.
 Applying Theorem
\ref{thm:rec1} to
$(\theta^*_i, \theta_i, i=0...d;  \varphi_j, \phi_{d-j+1}, j=1...d)$
and 
$f^*_0, f^*_1, \ldots, f^*_d$ we find that for $0 \leq i,j\leq d$,
\begin{eqnarray}
\theta^*_j f^*_i(\theta^*_j) = c^*_if^*_{i-1}(\theta^*_j)
+
a^*_i f^*_i(\theta^*_j)
+ b^*_i
f^*_{i+1}(\theta^*_j),
\label{eq:threetermd}
\end{eqnarray}
where $f^*_{-1}, f^*_{d+1}$ are indeterminates.
Comparing 
(\ref{eq:fipoly}) and
(\ref{eq:dualpoly1}) 
we find
\begin{eqnarray}
f_i(\theta_j)= f^*_j(\theta^*_i) \qquad \qquad  (0 \leq i,j\leq d).
\label{eq:dual}
\end{eqnarray}
Evaluating 
(\ref{eq:threetermd}) using 
(\ref{eq:dual}) and reindexing the result we obtain
(\ref{eq:dualrec}).
\hfill $\Box $ \\

\noindent
We finish this section  with a comment.

\begin{lemma}
\label{lem:kisalt}
With reference to Definition
\ref{def:ki} and 
Definition
\ref{def:bicidual},
\begin{eqnarray*}
k^*_i = \frac{
b^*_0b^*_1\cdots b^*_{i-1}
}
{
c^*_1c^*_2\cdots c^*_i
}
\qquad \qquad (0 \leq i \leq d).
\end{eqnarray*}
\end{lemma}
\noindent {\it Proof:}
Similar to the proof of Lemma 
\ref{lem:kialt}.
\hfill $\Box $ \\

\section{Some useful formulae}

In this section we give alternative  formulae for
the scalars $a_i, b_i, c_i$  from
Definition \ref{def:bici}. 
To avoid trivialities we assume the diameter
$d\geq 1$.
We begin with the $a_i$.

\begin{theorem} \cite[Lemma 5.1]{LS99}
\label{lem:aialt}
With reference to Definition \ref{def:bici}, let us assume $d\geq 1$.
Then
\begin{eqnarray}
a_0 &=& \theta_0+
\frac{\varphi_1}{\theta^*_0-\theta^*_1},
\\
a_i &=& \theta_i+
\frac{\varphi_i}{\theta^*_i-\theta^*_{i-1}}
+
\frac{\varphi_{i+1}}{\theta^*_i-\theta^*_{i+1}}
\qquad \qquad (1 \leq i \leq d-1),
\\
a_d &=& \theta_d+
\frac{\varphi_d}{\theta^*_d-\theta^*_{d-1}}.
\end{eqnarray}
\end{theorem}

\begin{lemma}
\label{lem:strec}
With reference to Definition \ref{def:bici}, assume $d\geq 1$.
Then
\begin{eqnarray}
c_i(\theta^*_{i-1}-\theta^*_i)-b_i(\theta^*_i-\theta^*_{i+1})= 
(\theta_1-\theta_0)(\theta^*_i-\theta^*_0)+\varphi_1
\label{eq:cibith}
\end{eqnarray}
for $0 \leq i \leq d$, where $\theta^*_{-1}, \theta^*_{d+1}$ 
denote indeterminates.
\end{lemma}
\noindent {\it Proof:}
Setting $\lambda =\theta_1$ in 
(\ref{eq:fipoly})
we find
\begin{eqnarray}
f_i(\theta_1) = 1 + 
\frac{(\theta_1-\theta_0)(\theta^*_i-\theta^*_0)}{\varphi_1}
\qquad \qquad (0 \leq i \leq d).
\label{eq:ths}
\end{eqnarray}
\noindent Setting $j=1$ in
(\ref{eq:threeterm}) and evaluating the 
result using 
(\ref{aiorig}), 
(\ref{eq:ths}) we  obtain 
(\ref{eq:cibith}).
\hfill $\Box $

\begin{theorem}
\label{lem:bicialt}
With reference to Definition \ref{def:bici}, assume $d\geq 1$.
Then
\begin{eqnarray}
b_0 
 &=& \frac
{
\varphi_1
}
{
\theta^*_{1}-\theta^*_0
},
\label{b0}
\\
b_i &=& 
\frac{
(\theta_0-a_i)(\theta^*_i-\theta^*_{i-1})+
(\theta_0-\theta_1)(\theta^*_0-\theta^*_i)+\varphi_1
}
{
\theta^*_{i+1}-\theta^*_{i-1}
}
\qquad  (1 \leq i \leq d-1),
\label{bi}
\\
b_d&=&0,
\label{bd}
\\
c_0&=&0, 
\label{c0}
\\
c_i &=& \frac
{
(\theta_0-a_i)(\theta^*_i-\theta^*_{i+1})+(\theta_0-\theta_1)(\theta^*_0-\theta^*_i)+\varphi_1
}
{
\theta^*_{i-1}-\theta^*_{i+1}
}
\qquad (1 \leq i \leq d-1),
\label{ci}
\\
c_d &=& \frac
{
\phi_d
}
{
\theta^*_{d-1}-\theta^*_d
}.
\label{cd}
\end{eqnarray}
\end{theorem}
\noindent {\it Proof:}
Lines (\ref{bd}), (\ref{c0}) are clear.
To get (\ref{b0})
set $i=0$ in
(\ref{biorig}).
To get (\ref{cd})
set $i=d$ in
(\ref{ciorig}).
To get (\ref{bi}), (\ref{ci}), solve the linear system
(\ref{aiorig}), 
(\ref{eq:cibith}) for $b_i, c_i$.
\hfill $\Box $
\\

\noindent Results similar to 
Theorem
\ref{lem:aialt},
Lemma \ref{lem:strec},
and Theorem 
\ref{lem:bicialt}
hold for the $a^*_i, b^*_i, c^*_i$.

\section{Remarks}

We conclude this paper with a few remarks.

\medskip
\noindent 
Let 
$(\theta_i, \theta^*_i, i=0...d;  \varphi_j, \phi_j, j=1...d)$
denote a parameter array over $\K$ and let
$f_0, f_1, \ldots, f_d$ denote the corresponding polynomials
from
Definition
\ref{def:assocpoly}.
Applying Theorem
\ref{eq:mth} with $\lambda=\theta_d $ 
and using
Lemma \ref{lem:scal}
we find
\begin{eqnarray}
f_i(\theta_d) = \frac{\phi_1 \phi_2\cdots \phi_i}
{\varphi_1 \varphi_2\cdots \varphi_i}
\qquad \qquad (0 \leq i \leq d).
\label{eq:fid}
\end{eqnarray}
Let the scalars $k_i$ be as in
Definition \ref{def:ki}.
Comparing 
(\ref{eq:fid})
with the formulae for $k_i$ given in 
Definition \ref{def:ki}
we find
\begin{eqnarray*}
k_i f_i(\theta_d) = 
\frac{
(\theta^*_0-\theta^*_1)
(\theta^*_0-\theta^*_2) \cdots 
(\theta^*_0-\theta^*_d)
}
{
(\theta^*_i-\theta^*_0) \cdots 
(\theta^*_i-\theta^*_{i-1}) 
(\theta^*_i-\theta^*_{i+1}) \cdots 
(\theta^*_i-\theta^*_d)
}
\end{eqnarray*}
for $0 \leq i \leq d$.

\medskip
\noindent 
We describe the matrix $G$ from
Theorem \ref{thm:gmat}. 
We use the following notation.
Let 
$(\theta_i, \theta^*_i, i=0...d;  \varphi_j, \phi_j, j=1...d)$
denote a parameter array over $\K$ and let
$q$ denote a base for this array.
To keep things simple we assume
$q \not=1$, $q\not=-1$.
For nonegative integers  $r,s,t$ such that $r+s+t\leq d$ we define
\begin{eqnarray*}
\lbrack r,s,t\rbrack_q:=
\frac{(q;q)_{r+s}
(q;q)_{r+t}
(q;q)_{s+t}}
{(q;q)_r (q;q)_s (q;q)_t (q;q)_{r+s+t}},
\end{eqnarray*}
where
\begin{eqnarray*}
(a;q)_n:=(1-a)(1-aq)\cdots (1-aq^{n-1}) \qquad \qquad n=0,1,2,\ldots
\end{eqnarray*}
We comment 
$\lbrack r,s,t\rbrack_q \in \K$
\cite[Definition 13.1]{LS24}.
Let $S$ denote the upper triangular
matrix in 
$\hbox{Mat}_{d+1}(\K)$
which has entries
\begin{eqnarray*}
S_{ij} = (\theta_0-\theta_d)
(\theta_0-\theta_{d-1})
\cdots (\theta_0-\theta_{d-j+i+1}) 
\lbrack i,j-i,d-j\rbrack_q
\end{eqnarray*}
for $0 \leq i\leq j\leq d$.
Then for  $G \in
 \mbox{Mat}_{d+1}(\K)$, 
$G$ satisfies
Theorem \ref{thm:gmat}(ii)
if and only if there exists a nonzero $\alpha \in \K$ such
that $G=\alpha S$
 \cite[Theorem 15.2]{LS24}.
 Similar results hold for $q = 1$ and
 $q=-1$ \cite[Lemma 13.2]{LS24}.

\section{Open problems}

\begin{problem}
Generalize 
Theorem \ref{eq:mth} so that it applies to polynomial sequences of infinite
length. Use  
this result to characterize the polynomials of the Askey scheme.
\end{problem}


\noindent Paul Terwilliger \hfil\break
Department of Mathematics \hfil\break
University of Wisconsin \hfil\break
480 Lincoln Drive \hfil\break
Madison, Wisconsin, 53706 USA 
\hfil\break
email: terwilli@math.wisc.edu \hfil\break

\end{document}